\documentclass[12pt,leqno]{article}

\usepackage{amsmath,amssymb,amsthm}
\usepackage[all]{xy}

\makeatletter

\swapnumbers

\newtheorem{theorem}{Theorem}[section]
\newtheorem{proposition}[theorem]{Proposition}
\newtheorem{lemma}[theorem]{Lemma}
\newtheorem{corollary}[theorem]{Corollary}

\theoremstyle{definition}

\newtheorem{definition}[theorem]{Definition}

\newtheorem{question}[theorem]{Question}

\theoremstyle{remark}

\newtheorem{remark}[theorem]{Remark}

\def\({{\rm (}}
\def\){{\rm )}}

\let\Mathrm\operator@font
\let\Cal\mathcal
\let\Bbb\mathbb
\newcommand{\fm}{\ensuremath{\mathfrak m}}
\newcommand{\fn}{\ensuremath{\mathfrak n}}
\newcommand{\fp}{\ensuremath{\mathfrak p}}

\def\standop#1{\mathop{\Mathrm #1}\nolimits}
\def\difstop#1#2{\expandafter\def\csname #1\endcsname{\standop{#2}}}
\def\defstop#1{\difstop{#1}{#1}}

\defstop{CMFI}
\defstop{Coker}
\defstop{Cl}

\defstop{Ext}

\difstop{height}{ht}
\defstop{Hom}

\defstop{Ker}

\defstop{Max}
\defstop{Min}
\defstop{MCM}

\defstop{NonSFR}
\defstop{NonCMFI}

\defstop{Reg}

\defstop{Spec}
\defstop{Supp}
\defstop{Sing}
\defstop{SFR}
\defstop{Soc}

\def\O{\Cal O}
\def\fm{\mathfrak{m}}
\def\fp{\mathfrak{p}}
\def\fq{\mathfrak{q}}
\def\fa{\mathfrak{a}}

\let\indlim\varinjlim
\let\projlim\varprojlim

\def\section{\@startsection{section}{1}{\z@ }%
{-3.5ex plus -1ex minus -.2ex}{2.3ex plus .2ex}{\bf }}

\long\def\refname{\par\kern -3ex
\begin{center}\rm R\sc{eferences}\end{center}\par\kern 
-2ex}

\def\@seccntformat#1{\csname the#1\endcsname.\quad}

\def\@@@sect#1#2#3#4#5#6[#7]#8{%
   \ifnum #2>\c@secnumdepth 
      \def \@svsec {}\else \refstepcounter {#1}%
      \def\@svsec{}
   \fi 
   \@tempskipa #5\relax 
   \ifdim \@tempskipa >\z@ 
     \begingroup #6\relax \@hangfrom {\hskip #3\relax 
     \@svsec}{\interlinepenalty \@M #8\par }\endgroup 
     \csname #1mark\endcsname {#7}
   \else 
   \def \@svsechd {#6\hskip #3\@svsec #8\csname #1mark\endcsname {#7}}
   \fi \@xsect {#5}}

\def\@@@startsection#1#2#3#4#5#6{%
 \if@noskipsec \leavevmode \fi \par \@tempskipa #4\relax \@afterindenttrue 
 \ifdim \@tempskipa <\z@ \@tempskipa -\@tempskipa \@afterindentfalse 
 \fi \if@nobreak \everypar {}\else \addpenalty {\@secpenalty }\addvspace 
  {\@tempskipa }\fi \@ifstar {\@ssect {#3}{#4}{#5}{#6}}{\@dblarg 
  {\@@@sect {#1}{#2}{#3}{#4}{#5}{#6}}}}

\def\theparagraph{\thesection.\arabic{paragraph}}
\def\aparagraph{\@@@startsection{paragraph}{2}{\z@ }%
              {1.75ex plus .2ex minus .15ex}{-1em}{\bf(\theparagraph) } }
\def\paragraph{\@@@startsection{paragraph}{2}{\z@ }%
              {1.75ex plus .2ex minus .15ex}{-1em}{}{\bf(\theparagraph)} }

\let\c@theorem\c@paragraph

\title{$F$-pure homomorphisms, strong $F$-regularity, and 
$F$-injectivity}

\author{M{\sc itsuyasu} H{\sc ashimoto}}
\date{\normalsize
Graduate School of Mathematics, Nagoya University\\
Chikusa-ku,  Nagoya 464--8602 JAPAN\\
{\small \tt hasimoto@math.nagoya-u.ac.jp}}

\begin{document}

\maketitle
\footnote[0]
    {Key words: $F$-pure, strongly $F$-regular, $F$-injective.
2010 \textit{Mathematics Subject Classification}. 
    Primary 13A35; Secondary 14L30.}

\begin{abstract}
We discuss Matijevic--Roberts type theorem on strong $F$-regularity,
$F$-purity, and Cohen--Macaulay $F$-injective (CMFI for short) property.
Related to this problem, we also discuss the 
base change problem and the openness of loci of these properties.
In particular, we define the notion of $F$-purity of homomorphisms
using Radu--Andr\'e homomorphisms, and prove basic properties of it.
We also discuss a strong version of strong $F$-regularity (very strong
$F$-regularity), and compare these two versions of strong $F$-regularity.
As a result, strong $F$-regularity and very strong $F$-regularity 
agree for local rings, $F$-finite rings,
and essentially finite-type algebras over an excellent local rings.
We prove the $F$-pure base change of strong $F$-regularity.
\end{abstract}

\section{Introduction}

Throughout this article, $p$ denotes a prime number.

The main objective of this paper is to prove Matijevic--Roberts type theorem
on strong $F$-regularity, $F$-purity, and Cohen--Macaulay $F$-injective
(CMFI, for short) properties.
$F$-purity was defined and discussed by 
M.~Hochster and J.~Roberts in 1970's \cite{HR}, \cite{HR2}.
It turned out that the $F$-purity is deeply connected with the notion of
log canonical singularities \cite{Watanabe}.
Strong $F$-regularity was defined by Hochster and Huneke \cite{HH3}
for $F$-finite rings.
$F$-injectivity was first defined by Fedder
\cite{Fedder}.
Recently, Schwede \cite{Schwede} proved that singularities of 
dense $F$-injective type in characteristic zero are Du Bois.

We prove 

\begin{trivlist}
\item[\bf Theorem~\ref{M-R.thm}]
Let $y$ be a point of $X$, and $Y$ the integral closed subscheme of $X$ whose
generic point is $y$.
Let $\eta$ be the generic point of an irreducible component of $Y^*$,
where $Y^*$ is the smallest $G$-stable closed subscheme of $X$ containing
$Y$.
Assume either that the second projection $p_2:G\times X\rightarrow X$ is
smooth, or that $S=\Spec k$ with $k$ a perfect field and $G$ is of finite
type over $S$.
Assume that $\O_{X,\eta}$ is of characteristic $p$.
Then $\O_{X,y}$ is of characteristic $p$.
Moreover, 
\begin{description}
\item[1] If $\O_{X,\eta}$ is $F$-pure, then $\O_{X,y}$ is $F$-pure.
\item[2] If $\O_{X,\eta}$ is excellent and strongly $F$-regular, then
$\O_{X,y}$ is strongly $F$-regular.
\item[3] If $\O_{X,\eta}$ is CMFI, then
$\O_{X,y}$ is CMFI.
\end{description}
\end{trivlist}

Matijevic--Roberts type theorems were originally conjectured and proved for
graded rings, see the introduction of \cite{HM} for a short history.
In \cite{HM}, roughly speaking, 
it is proved that if a ring theoretic property $\Bbb P$ enjoys
\lq smooth base change' and `flat descent,' then Matijevic--Roberts type
theorem for $\Bbb P$ under the action of a smooth group scheme holds.
Applying this principle, Matijevic--Roberts type theorems for
(weak) $F$-regularity and $F$-rationality were proved in \cite{HM}.

Smooth base change of $F$-purity is not so difficult.
In order to discuss this problem, 
we define the notion of $F$-purity of homomorphism
between (noetherian) commutative rings of characteristic $p$.
We use Radu--Andr\'e homomorphism to do so.
This map is used to characterize the regularity of 
homomorphism between
noetherian commutative rings of characteristic $p$.
Radu \cite{Radu} and Andr\'e \cite{Andre} proved that a 
homomorphism $f:A\rightarrow B$ between
noetherian commutative rings of characteristic $p$ is regular
(i.e., flat with geometrically regular fibers)
if and only if the Radu--Andr\'e homomorphism 
$\Phi_e(A,B):B^{(e)}\otimes_{A^{(e)}}A\rightarrow B$ 
given by $\Phi_e(A,B)(b^{(e)}\otimes a)=b^{p^e}a$ is
flat for some (or equivalently, any) $e>0$.
After that, the Radu--Andr\'e homomorphisms were used to study
the reducedness of homomorphisms by Dumitrescu, and
Cohen--Macaulay $F$-injective property by Enescu \cite{Enescu} and
the author \cite{Hashimoto2}.

In this article,
we define a homomorphism $f:A\rightarrow B$ between commutative rings
of characteristic $p$ is $F$-pure if $\Phi_e(A,B)$ is pure for any $e>0$.
This property behaves well under composition, localization, and base change,
and this notion is a canonical generalization of the notion of $F$-purity
of a ring.

Many properties of rings are promoted to properties of homomorphisms,
using the properties of (geometric) fibers.
For example, Grothendieck \cite[(6.8.1)]{EGA-IV} 
defined that a ring homomorphism $f:A\rightarrow B$
is Cohen--Macaulay if $f$ is flat with Cohen--Macaulay fibers.
However, it seems that it is not appropriate to define an $F$-pure
homomorphism to be a flat homomorphism with geometrically $F$-pure fibers,
because Singh \cite{Singh} constructed 
an example of homomorphism $f:A\rightarrow B$ of noetherian
rings of characteristic $p$ such that $A$ is a DVR (in particular, $F$-pure), 
$f$ is flat with 
geometrically $F$-pure fibers, but $B$ is not $F$-pure.

We also define the notion corresponding to the reduced property.
We say that a homomorphism $f:A\rightarrow B$ between commutative rings
of characteristic $p$ is Dumitrescu if $\Phi_e(A,B)$ is $A$-pure for any
$e>0$, see (\ref{Dumitrescu.par}).

Dumitrescu proved that if $f:A\rightarrow B$ is a flat homomorphism
between noetherian rings of characteristic $p$, then $f$ is Dumitrescu
if and only if $f$ is a reduced homomorphism 
(that is, a flat homomorphism with geometrically reduced fibers).

It is natural to ask whether a Dumitrescu homomorphism is flat.
We prove that if $A$ is regular, then a Dumitrescu homomorphism 
between noetherian commutative rings of characteristic $p$
$f:A\rightarrow B$ is flat.
A homomorphism $f:A\rightarrow B$ is said to be almost quasi-finite
if $f$ has finite fibers.
We prove that an almost quasi-finite Dumitrescu homomorphism is flat
(Theorem~\ref{Watanabe.thm}).

In the late 1980's, Hochster and Huneke defined $F$-regularity
using tight closure \cite{HH}.
They also defined strong $F$-regularity using Frobenius splittings
for $F$-finite rings of characteristic $p$ \cite{HH3}.
Hochster and Huneke defined the strong $F$-regularity for 
non-$F$-finite homomorphisms \cite[(5.3)]{HH2}.
Recently, Hochster \cite{Hochster} gave another definition of 
strong $F$-regularity.
In this paper, we call Hochster--Huneke's definition the very strong 
$F$-regularity (Definition~\ref{HH.def}), and Hochster's definition the
strong $F$-regularity.
We compare these two definitions.
Obviously, very strong $F$-regularity implies the strong $F$-regularity.
They
agree for local rings, $F$-finite rings, and essentially finite-type
algebras over excellent local rings.
We give a sufficient condition for the strong $F$-regular locus to be open
(Proposition~\ref{open.thm}).
We prove the $F$-pure base change of the strong $F$-regularity
(Theorem~\ref{Fedder-Watanabe.thm}).

We also discuss some basic properties of Cohen--Macaulay $F$-injectivity.
The base change of $F$-injectivity was first proved by
Aberbach--Enescu \cite{AE}, see \cite{Enescu2}.
We give another proof using Radu--Andr\'e homomorphism
(Proposition~\ref{base-change-CMFI.thm}).
This is a slight modification of Enescu's base change theorem on
$F$-rationality \cite{Enescu}.
We also prove the openness of CMFI locus (Corollary~\ref{CMFI-open.thm}).

In section~2, we discuss $F$-purity of homomorphisms.
In section~3, we prove basic properties of strong $F$-regularity and
very strong $F$-regularity, and prove the $F$-pure base change of
strong $F$-regularity.
In section~4, we discuss some properties of Cohen--Macaulay $F$-injectivity.
In section~5, we prove Matijevic--Roberts type theorem for $F$-purity,
strong $F$-regularity, and CMFI property.

The author thanks Professor K.-i. Watanabe for communicating the author
with his result (see Remark~\ref{Watanabe.rem}).
Special thanks are also due to 
Professor A. Singh and 
Professor K.-i. Yoshida for valuable advice.
The author is grateful to Professor K. Schwede and Professor F. Enescu 
for giving valuable comments to the former version of this paper.

\section{$F$-pure homomorphism}

\paragraph\label{Frobenius-notation.par}
Let $p$ be a prime number.
Let $k$ be a perfect field of characteristic $p$.
For a $k$-algebra $A$ of $k$ and $r\in\Bbb Z$, we define a new $k$-algebra
$A^{(r)}$ as follows.
$A^{(r)}$ as a ring is $A$, and the structure map of $A^{(r)}$ as a 
$k$-algebra is the composite
\[
k\xrightarrow{F^{-r}_k}k\xrightarrow{\varphi}A,
\]
where $\varphi$ is the original structure map of the $k$-algebra $A$,
and $F^{-r}_k(\alpha)=\alpha^{p^{-r}}$ for $\alpha\in k$.
The element $a\in A$, viewed as an element of $A^{(r)}$, is denoted by
$a^{(r)}$.
Note that for $e\geq 0$, the $e$th Frobenius map $F^e:A^{(e+r)}\rightarrow
A^{(r)}$ given by $F^e(a^{(e+r)})=(a^{p^e})^{(r)}$ is a $k$-algebra map.
This notation is consistent with that of Frobenius twisting 
in representation theory \cite[(I.9.10)]{Jantzen}.
Note that $A^{(-r)}$ is also written as ${}^r \!A$, $A^{(r)}$, 
$A^{(p^r)}$, or $A^{p^{-r}}$
by some authors.
For an $A$-module $M$, the module $M$ viewed as an $A^{(r)}$-module
(because $A^{(r)}=A$ as a ring) is denoted by $M^{(r)}$.
An element $m\in M$ viewed as an element of $M^{(r)}$ is denoted by
$m^{(r)}$.
If $e\geq 0$, then $M^{(-e)}$ is also an $A$-module through 
$F_A:A\rightarrow A^{(-e)}$.
The action of $A$ on $M^{(-e)}$ is given by 
$a\cdot m^{(-e)}=(a^{p^e}m)^{(-e)}$.

For a $k$-algebra map $f:A\rightarrow B$, the map $f$, viewed as a map
$A^{(r)}\rightarrow B^{(r)}$ is denoted by $f^{(r)}$.
Note that $f^{(r)}$ is a $k$-algebra map, and the diagram
\[
\xymatrix{
A^{(e)} \ar[r]^{f^{(e)}} \ar[d]^{F^e_A} & B^{(e)} \ar[d]^{F^e_B} \\
A \ar[r]^f & B
}
\]
is commutative for $e\geq 0$.
Thus, the $e$th Radu-Andr\'e homomorphism
\[
\Phi_e(A,B): B^{(e)}\otimes_{A^{(e)}}A\rightarrow B
\qquad
(b^{(e)}\otimes a \mapsto b^{p^e}(fa))
\]
is induced.
It is a $k$-algebra map for $e\geq 0$.
The $(-e)$-shift $\Phi_e(A,B)^{(-e)}: B\otimes_A A^{(-e)}\rightarrow B^{(-e)}
$ is denoted by $\Psi_e(A,B)$.

\begin{lemma}\label{pure-localization.thm}
Let $A\rightarrow B$ be a ring homomorphism,
and $\varphi:M\rightarrow M'$ a $B$-linear map between $B$-modules.
Then $\varphi$ is $B$-pure 
if and only if $\varphi_\fm:M_\fm\rightarrow M'_\fm$ is
$B_\fm$-pure for any maximal ideal $\fm$ of $A$.
\end{lemma}

\proof If $W$ is a $B_\fm$-module, then 
$1_W\otimes \varphi_\fm: W\otimes_{B_\fm}M_\fm\rightarrow W\otimes_{B_\fm}M'
_\fm$ is identified with $1_W\otimes \varphi: W\otimes_B M\rightarrow 
W\otimes_B M'$.
So the `only if' part is clear.

We prove the `if' part.
Let $N$ be any $B$-module, and let $K:=\Ker(1_N\otimes\varphi)$.
Then
\[
1_{A_\fm}\otimes (1_N\otimes\varphi):
A_\fm\otimes_A (N\otimes_B M)\rightarrow A_\fm\otimes_A(N\otimes_B M)
\]
is identified with
\[
1_{N_\fm}\otimes \varphi_\fm:
N_\fm\otimes_{B_\fm}M_\fm \rightarrow N_\fm\otimes_{B_\fm}M'_\fm.
\]
So $K_\fm=0$ for any $\fm$, and thus $K=0$ for any $N$.
This shows that $\varphi$ is $B$-pure.
\qed

\paragraph Consider that $k=\Bbb F_p$ in 
(\ref{Frobenius-notation.par}).
We say that a $k$-algebra map $f:A\rightarrow B$ is $F$-pure if
for any $e\geq 1$, $\Phi_e(A,B)$ is a pure ring homomorphism.

\begin{proposition}\label{basic.thm}
Let $f:A\rightarrow B$ and $g:B\rightarrow C$ be $\Bbb F_p$-algebra maps.
\begin{description}
\item[1] If $f$ and 
$g$ are $F$-pure homomorphisms, then so is $gf$.
\item[2] If $gf$ is $F$-pure and $g$ is pure, 
then $f$ is $F$-pure.
\item[3] An $\Bbb F_p$-algebra $A$ is $F$-pure if and only if the unique map
$\Bbb F_p\rightarrow A$ is $F$-pure.
\item[4] If $A$ is $F$-pure and $f$ is an $F$-pure homomorphism,
then $B$ is $F$-pure.
\item[5] A pure subring of an $F$-pure ring of characteristic $p$ is
again $F$-pure.
\item[6] A regular homomorphism of noetherian rings of characteristic $p$ 
is $F$-pure.
\item[7] A base change of an $F$-pure ring homomorphism is again $F$-pure.
Namely, if $f:A\rightarrow B$ is an $F$-pure homomorphism,
$A'$ an $A$-algebra, and $B':=B\otimes_A A'$, then $f\otimes1:A'\rightarrow
B'$ is $F$-pure.
\item[8] If 
$A\rightarrow A'$ is a pure ring homomorphism and 
$f\otimes 1:A'\rightarrow B':=B\otimes_AA'$ is $F$-pure,
then $f$ is $F$-pure.
\item[9] The following are equivalent.
\begin{description}
\item[(i)] $f$ is $F$-pure.
\item[(ii)] For any prime ideal $\fp$ of $A$, $f_\fp:A_\fp\rightarrow B_\fp$
is $F$-pure.
\item[(iii)] For any maximal ideal $\fm$ of $A$, $f_\fm:A_\fm\rightarrow B_\fm$
is $F$-pure.
\item[(iv)] For any prime ideal $\fq$ of $B$, $A\rightarrow B_\fq$ is
$F$-pure.
\item[(v)] For any maximal ideal $\fn$ of $B$, $A\rightarrow B_\fn$ is
$F$-pure.
\item[(vi)] For any prime ideal $\fq$ of $B$, $A_\fp\rightarrow B_\fq$ is
$F$-pure, where $\fp:=\fq\cap A$.
\item[(vii)] For any maximal ideal $\fn$ of $B$, $A_\fp\rightarrow B_\fn$ is
$F$-pure, where $\fp:=\fn\cap A$.
\end{description}
\end{description}
\end{proposition}

\proof {\bf 1} and {\bf 2} are 
obvious by \cite[Lemma~4.1, {\bf 1}]{Hashimoto}.
{\bf 3} follows from \cite[Lemma~4.1, {\bf 5}]{Hashimoto}.
{\bf 4} follows immediately by {\bf 1} and {\bf 3}.
{\bf 5} follows immediately by {\bf 2} and {\bf 3}.
{\bf 6} follows from a theorem of Radu \cite{Radu} and Andr\'e \cite{Andre}
which states that a homomorphism of noetherian rings of characteristic
$p$ $A\rightarrow B$ is regular if and only if there exists some $e>0$
such that $\Phi_e(A,B)$ is flat if and only if $\Phi_e(A,B)$ is faithfully
flat for any $e>0$, see also \cite{Dumitrescu}.
{\bf 7} and {\bf 8} follow 
easily from \cite[Lemma~4.1, {\bf 4}]{Hashimoto}.

{\bf 9} {\bf (ii)$\Rightarrow $(iii)}, {\bf (iv)$\Rightarrow$(v)},
and {\bf (vi)$\Rightarrow$(vii)} are trivial.
{\bf (i)$\Rightarrow$(ii)} is a consequence of {\bf 7}.
$\Phi_e(A,B)$ is pure if and only if $\Phi_e(A,B)_\fm$ is pure
for every $\fm\in\Max(A)$ by Lemma~\ref{pure-localization.thm},
where $\Max(A)$ denotes the set of maximal ideals of $A$.
Now {\bf (iii)$\Rightarrow$(i)} follows from 
\cite[Lemma~4.1, {\bf 4}]{Hashimoto}.
Applying \cite[Lemma~4.1, {\bf 1}]{Hashimoto} 
to $A\rightarrow B$ and $B
\rightarrow B_\fq$, we have that $\Phi_e(A,B)_{\fq^{(e)}}$ 
is pure if and only if
$\Phi_e(A,B_\fq)$ is pure, since $\Phi_e(B,B_\fq)$ is an isomorphism.
In view of Lemma~\ref{pure-localization.thm} again, 
this proves {\bf (i)$\Rightarrow$(iv)} and {\bf (v)$\Rightarrow$(i)}.
Combining {\bf (i)$\Rightarrow$(ii)} and {\bf (i)$\Rightarrow$(iv)}, 
we get {\bf (i)$\Rightarrow$(vi)} easily.
We prove {\bf (vii)$\Rightarrow$(v)}.
Applying \cite[Lemma~4.1, {\bf 1}]{Hashimoto} to $A\rightarrow A_\fp$
and $A_\fp\rightarrow B_\fn$, $\Phi_e(A,B_\fn)$ is pure if and only 
if $\Phi_e(A_\fp,B_\fn)$ is pure, since $\Phi_e(A,A_\fp)$ is an isomorphism.
The assertion is now clear.
\qed

\begin{lemma}\label{basic2.thm}
Let $f:A\rightarrow B$ be an $\Bbb F_p$-algebra map.
\begin{description}
\item[1] If $\Phi_1(A,B)$ is pure, then $f$ is $F$-pure.
\item[2] If $\Phi_e(A,B)$ is pure for some $e>0$ and $A$ is $F$-pure,
then $f$ is $F$-pure, and
$B^{(e')}\otimes_{A^{(e')}}A$ is $F$-pure for any $e'>0$.
\item[3] If $e>0$ and $B^{(e)}\otimes_{A^{(e)}}A$ is $F$-pure, then
$\Phi_e(A,B)$ is pure.
In particular, if $B^{(1)}\otimes_{A^{(1)}}A$ is $F$-pure, then $f$ is
$F$-pure.
\end{description}
\end{lemma}

\proof
{\bf 1} is an immediate consequence of \cite[Lemma~4.1, {\bf 2}]{Hashimoto}.
{\bf 2} $f$ is $F$-pure by {\bf 1} and \cite[Lemma~4.1, {\bf 2}]{Hashimoto}.
So $\Phi_{e'}(A,B)$ is $F$-pure, and $A$ is $F$-pure.
Now applying \cite[Lemma~4.1, {\bf 7}]{Hashimoto}, we have that
$B^{(e')}\otimes_{A^{(e')}}A$ is $F$-pure.
{\bf 3} By \cite[Lemma~4.1, {\bf 7}]{Hashimoto}, $\Phi_e(A,B)$ is pure.
The last assertion follows from {\bf 1}.
\qed

\begin{question}
Is a homomorphism $f:A\rightarrow B$ between noetherian rings of 
characteristic $p$ $F$-pure if $\Phi_e(A,B)$ is pure for some $e>0$?
\end{question}

\paragraph\label{Dumitrescu.par}
We say that a ring homomorphism $f:A\rightarrow B$ between
rings of characteristic $p$ is {\em Dumitrescu} if 
$\Phi_e(A,B)$ is pure as an $A$-linear map for every $e>0$.
By definition, an $F$-pure homomorphism is Dumitrescu.

If $\Phi_1(A,B)$ is $A$-pure, then $f$ is Dumitrescu by 
\cite[Lemma~4.1, {\bf 2}]{Hashimoto}.

Dumitrescu \cite{Dumitrescu2} proved that a {\em flat} ring homomorphism
$f:A\rightarrow B$ between noetherian rings of characteristic $p$ is
Dumitrescu if and only if $f$ is reduced, that is, $f$ is flat with 
geometrically reduced fibers.

\begin{proposition}\label{basic3.thm}
Let $f:A\rightarrow B$ and $g:B\rightarrow C$ be $\Bbb F_p$-algebra
maps.
\begin{description}
\item[1] If $f$ is $F$-pure and $g$ is Dumitrescu, then $gf$ is Dumitrescu.
\item[1'] If $f$ and $g$ are Dumitrescu, and $g$ is flat, then $gf$ is
Dumitrescu.
\item[2] If $gf$ is Dumitrescu and $g$ is $A$-pure, then $f$ is 
Dumitrescu.
\item[3] An $\Bbb F_p$-algebra $A$ is reduced if and only if the
unique map $\Bbb F_p\rightarrow A$ is Dumitrescu.
\item[4] If $A$ is $F$-pure and $f$ is Dumitrescu, then $B$ is reduced.
\item[5] A subring of a reduced ring \(of characteristic $p$\) is reduced.
\item[6] A reduced homomorphism of noetherian rings of characteristic $p$ 
is Dumitrescu.
\item[7] A base change of a Dumitrescu homomorphism is again Dumitrescu.
Namely, if $f:A\rightarrow B$ is a Dumitrescu homomorphism, $A'$ an
$A$-algebra, and $B':=B\otimes_AA'$, then $f\otimes 1:A'\rightarrow B'$
is Dumitrescu.
\item[8] If $A\rightarrow A'$ is a pure ring homomorphism and
$f\otimes 1:A'\rightarrow B':=B\otimes_AA'$ is Dumitrescu, then 
$f$ is Dumitrescu.
\item[9] The following are equivalent.
\begin{description}
\item[(i)] $f$ is Dumitrescu.
\item[(ii)] For any prime ideal $\fp$ of $A$, $f_\fp:A_\fp\rightarrow
B_\fp$ is Dumitrescu.
\item[(iii)] For any maximal ideal $\fm$ of $A$, $f_\fm:A_\fm
\rightarrow B_\fm$ is Dumitrescu.
\item[(iv)] For any prime ideal $\fq$ of $B$, $A\rightarrow B_\fq$ is
Dumitrescu.
\item[(v)] For any maximal ideal $\fn$ of $B$, $A\rightarrow B_\fn$
is Dumitrescu.
\item[(vi)] For any prime ideal $\fq$ of $B$, $A_\fp\rightarrow B_\fq$
is Dumitrescu, where $\fp:=\fq\cap A$.
\item[(vii)] For any maximal ideal $\fn$ of $B$, $A_\fp\rightarrow 
B_\fn$ is Dumitrescu, where $\fp:=\fn\cap A$.
\end{description}
\end{description}
\end{proposition}

\proof Similar to Proposition~\ref{basic.thm}.
\qed

\begin{lemma}\label{basic4.thm}
Let $f:A\rightarrow B$ be an $\Bbb F_p$-algebra map.
\begin{description}
\item[1] If $\Phi_1(A,B)$ is $A$-pure, then $f$ is Dumitrescu.
\item[2] If $\Phi_e(A,B)$ is $A$-pure for some $e>0$ and $A$ is $F$-pure,
then $f$ is Dumitrescu, and $B^{(e')}\otimes_{A^{(e')}}A$ is reduced for
any $e'>0$.
\end{description}
\end{lemma}

\proof Similar to Lemma~\ref{basic2.thm}.
\qed

\begin{question}
Is an $F$-pure homomorphism between noetherian rings of characteristic $p$
flat?
More generally, is a Dumitrescu homomorphism between noetherian rings
of characteristic $p$ flat?
\end{question}

\begin{lemma}\label{trivial.thm}
Let $f:A\rightarrow B$ be a ring homomorphism between rings of
characteristic $p$.
Then for $e\geq 0$, the composite
\[
A\cong A^{(e)}\otimes_{A^{(e)}}A
\xrightarrow{f^{(e)}\otimes 1}
B^{(e)}\otimes_{A^{(e)}}A
\xrightarrow{\Phi_e(A,B)}
B
\]
is $f$.
\end{lemma}

\proof Clear.
\qed

\begin{lemma}\label{Dumitrescu-pure.thm}
Let $f:A\rightarrow B$ be a ring homomorphism between rings of
characteristic $p$.
Assume that $A$ is noetherian, 
and
the image of the associated map
${}^af:\Spec B\rightarrow \Spec A$ contains $\Max(A)$, the set of 
maximal ideals of $A$.
If $f$ is Dumitrescu, then $f$ is pure.
\end{lemma}

\proof We may assume that $(A,\fm)$ is local.
Let $E$ be the injective hull of the residue field $A/\fm$ of $A$.
Set $E_n:=0:_E \fm^n$.
Let $A_n:=A/\fm^n$ and $B_n:=A_n\otimes_AB$.

It suffices to show that $f_n: A_n\rightarrow B_n$ is pure.
Indeed, if so, $E_n=E_n\otimes_AA\rightarrow E_n\otimes_AB$ is injective,
and hence taking the inductive limit, 
$E=E\otimes_AA\rightarrow E\otimes_AB$ is injective, and hence
$f$ is pure, see \cite{HR}.

So we may and shall assume that $(A,\fm)$ is artinian local.
Take $e>0$ sufficiently large so that $\fm^{[p^e]}=0$.
Namely, $a^{p^e}=0$ for every $a\in\fm$.
Then the Frobenius map $F^e:A^{(e)}\rightarrow A$ factors through
$(A/\fm)^{(e)}$.
Thus
\[
A=A^{(e)}\otimes_{A^{(e)}}A\xrightarrow{f^{(e)}\otimes 1}
B^{(e)}\otimes_{A^{(e)}}A
\]
is identified with the map
\[
A=(A/\fm)^{(e)}\otimes_{(A/\fm)^{(e)}}A
\xrightarrow{f_1^{(e)}\otimes 1}
(B/\fm B)^{(e)}\otimes_{(A/\fm)^{(e)}}A.
\]
Since $B/\fm B\neq 0$, this map is faithfully flat and hence is pure.

By the assumption and Lemma~\ref{trivial.thm}, $f$ is pure.
\qed

\begin{corollary}\label{local-pure.thm}
A local homomorphism between noetherian local rings of characteristic
$p$ is pure, if it is Dumitrescu.
\end{corollary}

\begin{lemma}\label{noetherian.thm}
If $f:A\rightarrow B$ is $F$-pure and $B$ is noetherian,
then $B^{(e)}\otimes_{A^{(e)}}A$ is noetherian for $e\geq 0$.
\end{lemma}

\proof This is because $\Phi_e(A,B)$ is pure.
\qed

\begin{lemma}
Let $K$ be a field of characteristic $p$, and $B$ a $K$-algebra.
Then the following are equivalent.
\begin{description}
\item[1] $K\rightarrow B$ is $F$-pure, and $B$ is noetherian.
\item[2] For any $e>0$, $B\otimes_K K^{(-e)}$ is noetherian and $F$-pure.
\item[3] There exists some $e>0$ such that $B\otimes_K K^{(-e)}$ is
noetherian and $F$-pure.
\item[4] $B$ is noetherian, and 
$B$ is geometrically $F$-pure over $K$, that is to say, for any finite 
algebraic extension $L$ of $K$, $B\otimes_K L$ is $F$-pure.
\end{description}
\end{lemma}

\proof Note that $(B\otimes_K K^{(-e)})^{(e)}\cong B^{(e)}\otimes_{K^{(e)}}K$.
{\bf 1$\Rightarrow$2}
Let $e>0$.
$B\otimes_K K^{(-e)}$ is noetherian by Lemma~\ref{noetherian.thm}.
$B\otimes_K K^{(-e)}$ is $F$-pure by Lemma~\ref{basic2.thm}, {\bf 2}.

{\bf 2$\Rightarrow $3} is trivial.

{\bf 3$\Rightarrow$1} 
As
\[
B=B\otimes_KK\xrightarrow{1_B\otimes F^e}
B\otimes_K K^{(-e)}
\]
is faithfully flat and $B\otimes_KK^{(-e)}$ is noetherian, 
$B$ is noetherian.
$K\rightarrow B$ is $F$-pure by Lemma~\ref{basic2.thm}, {\bf 3}.

{\bf 1,2,3$\Rightarrow$4} $B$ is noetherian, as assumed.
As the field $L$ is $F$-pure and $L\rightarrow B\otimes_K L$ is
$F$-pure (as it is the base change of the $F$-pure homomorphism
$K\rightarrow B$), $B\otimes_K L$ is also $F$-pure by
Proposition~\ref{basic.thm}, {\bf 4}.

{\bf 4$\Rightarrow$ 1} Let $L$ be a finite extension field of $K$ such that 
$L\subset K^{(-1)}$.
Then $F:B^{(1)}\otimes_{K^{(1)}}L^{(1)}\rightarrow B\otimes_K L$ is pure.
As $F$ factors through $B$, the map
$B^{(1)}\otimes_{K^{(1)}}L^{(1)}\rightarrow B$ is pure.
Taking the inductive limit on $L$, 
$\Phi_1(K,B):B^{(1)}\otimes_{K^{(1)}}K\rightarrow B$ is also pure.
So $K\rightarrow B$ is $F$-pure.
\qed

\begin{corollary}
If $f:A\rightarrow B$ is $F$-pure homomorphism between noetherian rings
of characteristic $p$, then $f$ has geometrically $F$-pure fibers.
That is, for any $\fp\in \Spec A$ and any finite algebraic extension
$L$ of $\kappa(\fp)$, $B\otimes_A L$ is $F$-pure.
\end{corollary}

\begin{remark}
The converse is not true in general.
Indeed, a flat homomorphism with geometrically $F$-pure fibers need
not be an $F$-pure homomorphism.
Singh's example \cite[section~6]{Singh} 
shows that for $p>2$, there is an example of a homomorphism $f:A\rightarrow B$
such that $A=\Bbb F_p[t]_{(t)}$, $f$ is flat with geometrically
$F$-pure fibers, but $B$ is not $F$-pure.
If $f$ were $F$-pure, then $B$ must have been $F$-pure by
Proposition~\ref{basic.thm}, {\bf 4}.
\end{remark}

\paragraph 
Let $f:A\rightarrow B$ be a homomorphism between noetherian rings.
We say that $f$ is almost quasi-finite if for any $P\in\Spec A$, 
$\kappa(P)\otimes_A B$ is module finite over $\kappa(P)$.
A quasi-finite homomorphism (that is, finite-type homomorphism with
zero-dimensional fibers) is almost quasi-finite.
A localization $A\rightarrow A_S$ is almost quasi-finite.
A composite of almost quasi-finite homomorphisms is almost quasi-finite.
A base change of an almost quasi-finite homomorphism is almost quasi-finite.

\begin{theorem}\label{Watanabe.thm}
Let $f:A\rightarrow B$ be an almost quasi-finite 
homomorphism between noetherian rings of characteristic $p$.
Then the following are equivalent.
\begin{description}
\item[1] $f$ is regular.
\item[2] $f$ is $F$-pure.
\item[3] $f$ is Dumitrescu.
\end{description}
\end{theorem}

\proof Note that {\bf 1$\Rightarrow$2$\Rightarrow$3} is
trivial.
So it suffices to prove {\bf 3$\Rightarrow$1}.
Note that each fiber $\kappa(P)\otimes_AB$ is Dumitrescu over $\kappa(P)$
by Proposition~\ref{basic3.thm}, {\bf 7}.
So it is geometrically reduced by Dumitrescu's theorem \cite{Dumitrescu2}.
As we assume that $\kappa(P)\otimes_A B$ is finite over $\kappa(P)$,
we have that $\kappa(P)\otimes_A B$ is \'etale over $\kappa(P)$.
So in order to prove {\bf 1}, it suffices to prove that $f$ is flat.

Thanks to Proposition~\ref{basic3.thm}, {\bf 9},
we may assume that $f:(A,\fm)\rightarrow (B,\fn)$ is a local homomorphism
between local rings.
By the local criterion of flatness 
((5)$\Rightarrow$(1) of \cite[Theorem~22.3]{CRT})
and Proposition~\ref{basic3.thm}, {\bf 7}, we may assume that
$A$ is artinian.
Then $B$ is module finite over $A$.
It is easy to see that
if $l_A(B)\geq l_A(A)l_A(B/\fm B)$, then $f$ is flat
by (4)$\Rightarrow$(1) of \cite[Theorem~22.3]{CRT},
where $l_A$ denotes the length as an $A$-module.

Take $e>0$ sufficiently large so that $\fm^{(e)}A=0$, that is, 
for any $x\in \fm$, $x^{p^e}=0$.
Then $\Phi_e(A,B)$ is identified with
\[
(B/\fm B)^{(e)}\otimes_{(A/\fm)^{(e)}}A\rightarrow B.
\]
This map is an injective $A$-linear map with
$l_A((B/\fm B)^{(e)}\otimes_{(A/\fm)^{(e)}}A)=l_A(A)l_A(B/\fm B)$.
So $l_A(B)\geq l_A(A)l_A(B/\fm B)$, as desired.
\qed

\begin{remark}\label{Watanabe.rem}
Theorem~\ref{Watanabe.thm} for the case that $f:A\rightarrow B$ is a
finite homomorphism between integral domains (the crucial case) is due to 
K.-i.~Watanabe.
\end{remark}

\begin{proposition}\label{Dumitrescu-flat.thm}
Let $f:(A,\fm)\rightarrow (B,\fn)$ be a local homomorphism between
noetherian local rings of characteristic $p$, and $t\in\fm$.
Assume that $A$ is normally flat along $tA$,
$A/tA$ is reduced, $A/tA\rightarrow B/tB$ is flat, and
$f$ is Dumitrescu.
Then $f$ is flat.
\end{proposition}

\proof It suffices to show that $A/t^nA\rightarrow B/t^nB$ is flat for
all $n\geq 1$.
So we may assume that $t^n=0$ for some $n\geq 1$.
We prove the assertion by induction on $n$.
If $n=1$, then the assertion is assumed by the assumption of the proposition.

So we consider the case that $n\geq 2$.
It suffices to show that the canonical map
\[
\gamma_i:B/tB\otimes_{A/tA}(t^iA/t^{i+1}A)\rightarrow t^iB/t^{i+1}B
\qquad
\bar b\otimes \overline{t^ia}\mapsto \overline{t^iab}
\]
is injective for $i=1,\ldots,n-1$.
This is true for $i=1,\ldots,n-2$, as $A/t^{n-1}A\rightarrow B/t^{n-1}B$ 
is flat by induction assumption.
So it suffices to show that $\gamma_{n-1}:B/tB\otimes_{A/tA}t^{n-1}A
\rightarrow t^{n-1}B$ is injective.
Let $x\in B$ be an element such that
$\bar x\otimes t^{n-1}\in\Ker\gamma_{n-1}$,
or equivalently, $xt^{n-1}=0$.
Take $e\gg0$ such that $p^e>n$.
Obviously, we have $x^{p^e}t^{n-1}=0$ in $B$.
As
\[
\Phi_e(A,B):(B/tB)^{(e)}\otimes_{(A/tA)^{(e)}}A\cong
B^{(e)}\otimes_{A^{(e)}}A\rightarrow B
\]
is injective, $\bar x^{(e)}\otimes t^{n-1}=0$ in $(B/tB)^{(e)}\otimes
_{(A/tA)^{(e)}}A$.
As $(B/tB)^{(e)}$ is $(A/tA)^{(e)}$-flat, $\bar x^{(e)}\otimes t^{n-1}=0$ in
$(B/tB)^{(e)}\otimes_{(A/tA)^{(e)}}t^{n-1}A$.
Since $t^{n-1}A$ is a rank-one free $A/tA$-module with $t^{n-1}$ its 
basis, $\bar 
x^{(e)}\otimes 1\in (B/tB)^{(e)}\otimes_{(A/tA)^{(e)}}A/tA$ is zero.
As $A/tA$ is reduced, $(A/tA)^{(e)}\rightarrow A/tA$ is injective.
So $(B/tB)^{(e)}\rightarrow 
(B/tB)^{(e)}\otimes_{(A/tA)^{(e)}}A/tA$ is injective.
Hence $\bar x^{(e)}=0$ in $(B/tB)^{(e)}$.
So $\bar x=0$ in $B/tB$.
This shows that $\bar x\otimes t^{n-1}=0$ in $B/tB\otimes_{A/tA}t^{n-1}A$,
and thus $\gamma_{n-1}$ is injective, as desired.
\qed

\begin{corollary}\label{Dumitrescu-flat2.thm}
Let $f:(A,\fm)\rightarrow (B,\fn)$ be a Dumitrescu local homomorphism
between noetherian local rings of characteristic $p$.
If $t\in \fm$ is a nonzerodivisor, $A/tA$ is reduced, 
and $B/tB$ is $A/tA$-flat, then $B$ is $A$-flat.
\qed
\end{corollary}

\begin{corollary}
Let $f:(A,\fm)\rightarrow (B,\fn)$ be a Dumitrescu local homomorphism
between noetherian local rings of characteristic $p$.
If $A$ is regular, then $f$ is flat.
\end{corollary}

\proof We prove this by induction on $\dim A$.
If $\dim A=0$, then $A$ is a field, and $f$ is flat.

Next consider the case $\dim A>0$.
Then take $t\in \fm\setminus \fm^2$.
Then $t$ is a nonzerodivisor, 
$A/tA$ is regular, and $A/tA\rightarrow B/tB$ is flat by induction 
assumption.
By Corollary~\ref{Dumitrescu-flat2.thm}, $f$ is flat.
\qed

\section{Strong $F$-regularity}

\paragraph For a ring $R$, we define 
$R^\circ:=R\setminus\bigcup_{P\in\Min R}P$, where $\Min R$ denotes the
set of minimal primes of $R$.
Let $R$ be a ring of characteristic $p$.
Let $M$ be an $R$-module, and $N$ its submodule.
We define
\begin{multline*}
\Cl_R(N,M)=N_M^*:=\{x\in M \mid \exists c\in R^\circ\;\exists e_0\geq 1\;\\
\forall e\geq e_0\; x\otimes c^{(-e)}\in M/N\otimes_RR^{(-e)} 
\text{ is zero}\},
\end{multline*}
and call it the tight closure of $N$ in $M$.
Note that $\Cl_R(N,M)$ is an $R$-submodule of $M$ containing $N$
\cite[section~8]{HH}.

\begin{lemma}\label{tight-closure.thm}
Let $R$ be a noetherian commutative ring of characteristic
$p$, and $S$ a multiplicatively closed subset of $R$.
Let $M$ be an $R_S$-module, and $N$ its $R_S$-submodule.
Then $\Cl_R(N,M)=\Cl_{R_S}(N,M)$.
\end{lemma}

\proof
Note that 
\begin{multline*}
x\in \Cl_R(N,M) \iff
\exists c\in R^\circ\;
\exists q'\;
\forall q\geq q'\;
x\otimes c^{(-e)}\in 
\\
\Ker(M\otimes_R R^{(-e)}\rightarrow 
M/N\otimes_R R^{(-e)})=:N^{[q]}
\end{multline*}
and that
\begin{multline*}
x\in \Cl_{R_S}(N,M) \iff
\exists c\in R_S^\circ\;
\exists q'\;
\forall q\geq q'\;
x\otimes c^{(-e)}\in 
\\
\Ker(M\otimes_{R_S} R_{S}^{(-e)}\rightarrow 
M/N\otimes_{R_S} R_S^{(-e)})=N^{[q]},
\end{multline*}
where $q=p^e$ and $q'$ denote some power of $p$.
If $c\in R^\circ$, then $c/1\in R_S^\circ$.
So $\Cl_R(N,M)\subset \Cl_{R_S}(N,M)$.

Let $c\in R$, $c/1\in R_S^\circ$, and assume that there exists some
$q'$ such that for all $q\geq q'$, $x\otimes c^{(-e)}\in N^{[q]}$.
Take $\delta\in R$ such that for any $P\in\Min R$, $\delta\in P$
if and only if $c\notin P$.
Then $\delta/1$ is nilpotent in $R_S$.
Replacing $\delta$ by its some power, we may assume that $\delta/1=0$
in $R_S$.
Then $x\otimes(c+\delta)^{(-e)}\in N^{[q]}$ 
for $q\geq q'$
and $c+\delta\in R^\circ$.
Hence $x\in \Cl_R(N,M)$.
This shows $\Cl_R(N,M)\supset \Cl_{R_S}(N,M)$.
\qed

\def\citinfo{Hochster \cite[p.~166]{Hochster}}
\begin{definition}[\citinfo]\label{Hochster.def}
We say that a noetherian ring $R$ of characteristic $p$ is
{\em strongly $F$-regular} if $\Cl_R(N,M)=N$ for any $R$-module $M$
and any submodule $N$ of $M$.
\end{definition}

\def\citinfo{cf.\ \cite[(5.3)]{HH2}}
\begin{definition}[\citinfo]\label{HH.def}
We say that a noetherian ring $R$ of characteristic $p$ is {\em very
strongly $F$-regular} if for any $c\in R^\circ$
there exists some $e>0$ such that the map
$c^{(-e)}F^e:R\rightarrow R^{(-e)}$ ($x\mapsto (cx^{p^e})^{(-e)}$)
is $R$-pure.
\end{definition}

\begin{lemma}\label{hr.thm}
Let $(R,\fm)$ be a local ring, and $E$ be the injective
hull of the residue field.
Let $S$ be an $R$-module, and $h:R\rightarrow S$ be an $R$-linear map.
If $1_E\otimes h:E\rightarrow E\otimes_R S$ is injective, then $h$ is
$R$-pure.
\end{lemma}

\proof Exactly the same proof as \cite[(6.11)]{HR} works.
\qed

\begin{lemma}\label{strongly-F-regular.thm}
Let $R$ be a noetherian ring of characteristic $p$.
Then the following are equivalent.
\begin{description}
\item[1] $R$ is strongly $F$-regular.
\item[2] For any multiplicatively closed subset $S$ of $R$, $R_S$ is
strongly $F$-regular.
\item[3] $R_\fm$ is strongly $F$-regular for any $\fm\in\Max(R)$.
\item[4] For any $\fm\in\Max(R)$, $\Cl_R(0,E_{R}(R/\fm))=0$.
\item[5] For any $\fm\in\Max(R)$, $\Cl_{R_\fm}(0,E_{R_\fm}(R/\fm))=0$.
\item[6] $R_\fm$ is very strongly $F$-regular for any $\fm\in\Max(R)$.
\end{description}
\end{lemma}

\proof {\bf 1$\Rightarrow$4}, {\bf 2$\Rightarrow$3}, 
and {\bf 3$\Rightarrow$5} are obvious.
{\bf 1$\Rightarrow$2} and {\bf 4$\Rightarrow$5} follow from 
Lemma~\ref{tight-closure.thm}.

We prove {\bf 5$\Rightarrow$6}.
We may assume that $(R,\fm)$ is local.
Let $E=E_R(R/\fm)$ be the injective hull of the residue field.
The kernel of the map
\[
1_E\otimes F_R: E\rightarrow E\otimes_R R^{(-1)}
\qquad
(x\mapsto x\otimes 1)
\]
is contained in $\Cl_R(0,E)=0$.
Thus $F_R:R\rightarrow R^{(-1)}$ is pure by \cite[(6.11)]{HR}.
In other words, $R$ is $F$-pure.

Now let $c\in R^\circ$, and set $K_e:=\Ker(1_E\otimes c^{(-e)}F^e:
E\rightarrow E\otimes_R R^{(-e)})$,
where $1_E\otimes c^{(-e)}F^e$ sends $x$ to $x\otimes c^{(-e)}$.
For $e'>e$, we have a commutative diagram
\begin{equation}\label{E.eq}
\xymatrix{
E~~ \ar[r]^-{1_E\otimes c^{(-e)}F^e} \ar[d]^{1_E\otimes c^{(-e')}F^{e'}} &
~~E\otimes_R R^{(-e)} \ar[d]^{1\otimes F^{e'-e}} \\
E\otimes_R R^{(-e')}~~~~~~ \ar[r]^{1_E\otimes(c^{p^{e'-e}-1})^{(-e')}} &
~~~~~~E\otimes_R R^{(-e')}
}.
\end{equation}
As $1\otimes F^{e'-e}$ is injective, we have 
\[
E\supset K_1\supset K_2\supset \cdots.
\]
As $E$ is an artinian module, there exists some $e\gg0$ such that
$K_e=\bigcap_{e'}K_{e'}\subset \Cl_R(0,E)=0$.
By Lemma~\ref{hr.thm}, $c^{(-e)}F^e:R\rightarrow R^{(-e)}$ is pure,
as desired.

{\bf 6$\Rightarrow$3} 
We may assume that $(R,\fm)$ is local.
Considering the case that $c=1$, $F^e:R\rightarrow R^{(-e)}$ is pure
for some $e$.
Hence $R$ is $F$-pure.
Now let $c\in R^\circ$.
Then there exists some $e>0$ such that $c^{(-e)}F^e$ is $R$-pure.
Considering the commutative diagram (\ref{E.eq}), we have that
$c^{(-e')}F^{e'}$ is $R$-pure for $e'\geq e$.

Now let $M$ be an $R$-module, $N$ its submodule.
Then
$1_{M/N}\otimes c^{(-e')}F^{e'}:M/N\rightarrow M/N\otimes_R R^{(-e')}$ 
is injective by the purity for $e'\geq e$.
This shows that $\Cl_R(N,M)=N$.

{\bf 3$\Rightarrow$1}
Let $M$ be an $R$-module and $N$ its submodule.
Let $\fm\in\Max(R)$.
Then
\[
N_\fm\subset \Cl_R(N,M)_\fm\subset \Cl_{R_\fm}(N_\fm,M_\fm)=N_\fm.
\]
Hence $\Cl_R(N,M)_\fm=N_\fm$ for any $\fm\in\Max(R)$.
This shows that $\Cl_R(N,M)=N$.
\qed

\begin{corollary}\label{sfr-f-pure.thm}
A strongly $F$-regular noetherian ring of characteristic $p$ is 
$F$-regular.
In particular, it is normal and $F$-pure.
\end{corollary}

\proof By the definition of strong $F$-regularity, a strongly $F$-regular
implies weakly $F$-regular.
By {\bf 1$\Rightarrow$2} of the lemma, it is also $F$-regular.

The normality assertion is a consequence of \cite[(4.2)]{HH2}.

For the $F$-purity assertion, it suffices to point out that a weakly
$F$-regular noetherian ring $R$ of characteristic $p$ is $F$-pure 
\cite[Remark~1.6]{FW}.
Almost by the definition of weak $F$-regularity, for any ideal $I$ of $R$,
$1\otimes F: R/I\otimes_R R\rightarrow R/I\otimes_R R^{(-1)}$ is
injective.
Thus $F_R:R\rightarrow R^{(-1)}$ is cyclically pure.
But as $R$ is normal, it is approximately Gorenstein, and thus $F_R$ is
pure \cite{Hochster2}.
\qed

\begin{lemma}\label{vsfr-localization.thm}
If $R$ is very strongly $F$-regular noetherian ring
of characteristic $p$, then $R_S$ is very strongly $F$-regular for any
multiplicatively closed subset $S$ of $R$.
In particular, a very strongly $F$-regular noetherian ring of
characteristic $p$ is strongly $F$-regular.
\end{lemma}

\proof Let $c/s\in (R_S)^\circ$, where $c\in R$ and $s\in S$.
Take $\delta\in R$ such that $\delta \in P$ if and only if $c\notin P$ for
$P\in \Min R$.
Replacing $\delta$ by its power, we may assume that $\delta/1=0$
in $R_S$.
Set $d=c+\delta$.
Then $d^{(-e)}F^e:R\rightarrow R^{(-e)}$ is $R$-pure for some $e\geq 0$.
So $d^{(-e)}F^e:R_S\rightarrow R_S^{(-e)}$ is $R_S$-pure.
It follows that 
$(c/s)^{(-e)}F^e:R_S\rightarrow R_S^{(-e)}$ is also $R_S$-pure,
since $c/s=d/s$.
The last assertion follows from Lemma~\ref{strongly-F-regular.thm}.
\qed

\begin{lemma}
Let $R$ be an $F$-finite noetherian ring of characteristic $p$.
If $R$ is strongly $F$-regular, then it is very strongly $F$-regular.
\end{lemma}

\proof $R_\fm$ is very strongly $F$-regular for $\fm\in\Max R$ by
Lemma~\ref{strongly-F-regular.thm}.
Then by \cite[(5.2)]{HR2}, $R_\fm$ is ``strongly $F$-regular''
in the sense of \cite{HH3} for each $\fm$.
By \cite[(3.1)]{HH3},
$R$ is ``strongly $F$-regular'' in the sense of \cite{HH3}.
So $R$ is very strongly $F$-regular.
\qed

\begin{lemma}\label{direct-product.thm}
Let $R=R_1\times R_2$ be a noetherian ring of characteristic
$p$.
Then $R$ is very strongly $F$-regular if and only if both $R_1$ and $R_2$ are.
\end{lemma}

\proof If $R$ is very strongly $F$-regular, then $R_1$ and $R_2$ are
very strongly $F$-regular, as $R_1$ and $R_2$ are localizations of $R$.

Conversely, let $R_1$ and $R_2$ be very strongly $F$-regular.
Let $e_1$ and $e_2$ be respectively the idempotents corresponding to $R_1$
and $R_2$.
Let $c\in R^\circ$.
Then $ce_i\in R_i^\circ$, and hence there exist some $r_i$ ($i=1,2$)
such that $(ce_i)^{(-r_i)}F^{r_i}:R_i\rightarrow R_i^{(-r_i)}$ is
$R_i$-pure for each $i$.
As each $R_i$ is $F$-pure, it is easy to see that letting $r=\max(r_1,r_2)$,
$(ce_i)^{(-r)}F^r:R_i\rightarrow R_i^{(-r)}$ is $R_i$-pure for $i=1,2$.
Then $c^{(-r)}F^r:R\rightarrow R^{(-r)}$ is $R$-pure.
\qed

\begin{lemma}\label{normality.thm}
Let $R\rightarrow S$ be a ring homomorphism.
Let $S=S_1\times\cdots\times S_r$ be a finite direct product of
integrally closed domains.
Assume that for any nonzerodivisor $a$ of $R$, $aS\cap R=aR$.
Then $R$ is integrally closed in the total quotient ring $Q(R)$.
\end{lemma}

\proof Let $\alpha=b/a$ be an element of $Q(R)$, where $a,b\in R$ with
$a$ a nonzero divisor of $R$.
Assume that $\alpha$ is integral over $R$.

Let $1\leq i\leq r$.
Consider the case that $a$ is nonzero in $S_i$.
Then $b/a$ makes sense in 
the field of fractions $Q(S_i)$ of $S_i$, 
and it is integral over $S_i$.
As $S_i$ is integrally closed, $b/a\in S_i$.
Hence $b\in aS_i$.

Next consider the case that $a$ is zero in $S_i$.
Since $b/a$ is integral over $R$ and $a$ is a nonzerodivisor, 
there exists some $n\geq 1$ such that $b^n\in aR$.
This shows that $b^n=0$ in $S_i$.
As $S_i$ is a domain, $b=0$ in $S_i$.
This shows that $b\in aS_i$.

As $b\in aS_i$ for any $i$, $b\in aS\cap R=aR$.
Hence $\alpha=b/a\in R$, and $R$ is integrally closed in $Q(R)$.
\qed

\begin{corollary}\label{cyclically-pure-normal.thm}
Let $S$ be a noetherian normal ring, and $R$ its cyclically pure
subring.
Then $R$ is a noetherian normal ring, and hence $R$ is a pure subring of $S$.
\end{corollary}

\proof If $I_1\subset I_2\subset\cdots$ is an ascending chain of ideals in
$R$, then $I_1S\subset I_2S\subset\cdots$ is that of $S$, and hence
$I_NS=I_{N+1}S=\cdots$ for  some $N$.
Hence $I_NS\cap R=I_{N+1}S\cap R=\cdots$.
By cyclic purity, $I_N=I_{N+1}=\cdots$, and hence  $R$ is noetherian.

As $R$ is a subring of $S$, $R$ is reduced.
By Lemma~\ref{normality.thm}, $R$ is integrally closed in $Q(R)$.
So $R$ is a noetherian normal ring.

Hence $R$ is approximately Gorenstein \cite{Hochster2}, and hence
$R$ is a pure subring of $S$.
\qed

\begin{corollary}\label{cyclically-pure-F-reg.thm}
Let $S$ be a noetherian ring of characteristic $p$, and
$R$ is cyclically pure subring of $S$.
If $S$ is weakly $F$-regular \(resp.\ $F$-regular\), then
so is $R$.
\end{corollary}

\proof In any case, $S$ is normal.
So $R$ is normal and pure in $S$.
Now the assertion follows from \cite[(9.11)]{HO}.
\qed

\begin{lemma}\label{pure-very-strongly-F-regular.thm}
Let $R\rightarrow S$ be a cyclically pure ring homomorphism
between noetherian rings of characteristic $p$.
If $S$ is very strongly $F$-regular, then so is $R$.
\end{lemma}

\proof Note that $S$ is $F$-regular.
So $R$ is normal and is a pure subring of $S$ by
Corollary~\ref{cyclically-pure-normal.thm}.
By Lemma~\ref{direct-product.thm}, we may assume that $R$ is a normal domain.
We can write $S=S_1\times\cdots\times S_r$, where $S_i$ is a very
strongly $F$-regular domain for each $i$.
Assume that the induced map $f_i:R\rightarrow S_i$ is not injective.
Let $I_i=\Ker f_i$.
For $\fm\in \Max R$, $(I_i)_\fm\neq 0$.
So $R_\fm\rightarrow (S_i)_\fn$ cannot be injective for any $\fn\in\Max S_i$.
By \cite[(9.10)]{HO}, $R_\fm\rightarrow S'_\fn$ is pure for some
maximal ideal $\fn$ of $S':=S_1\times\cdots\times S_{i-1}\times S_{i+1}
\times\cdots\times S_r$.
This shows that $R\rightarrow S'$ is still pure.
Removing redundant $S_i$, we may assume that $f_i$ is injective for each
$i$.

Now let $c\in R^\circ$.
Then by our additional assumption, $c\in S^\circ$.
So $c^{(-e)}F^e_S$ is $S$-pure for some $e\geq 1$.
As the diagram
\[
\xymatrix{
~~S~~ \ar[r]^-{c^{(-e)}F^e_S} & ~~S^{(-e)} ~~\\
R\ar[r]^-{c^{(-e)}F^e_R} \ar[u] & R^{(-e)} \ar[u]
}
\]
is commutative, we have that $c^{(-e)}F^e_R$ is $R$-pure.
Hence $R$ is very strongly $F$-regular.
\qed

\begin{lemma}\label{affine-open.thm}
Let $R\rightarrow R'$ be a homomorphism between 
noetherian rings of characteristic $p$.
Assume that the induced map $\Spec R'\rightarrow \Spec R$ is an open immersion.
If $R$ is very strongly $F$-regular, then $R'$ is very strongly $F$-regular.
\end{lemma}

\proof $\Spec R'$ has a finite affine open covering $\Spec R'=\bigcup U_i$
such that $U_i=\Spec R[1/f_i]$ for some $f_i\in R$.
Then $R[1/f_i]$ is very strongly $F$-regular by
Lemma~\ref{vsfr-localization.thm}.
So $R'':=\prod_i R[1/f_i]$ is also very strongly $F$-regular by
Lemma~\ref{direct-product.thm}.
As $R''$ is faithfully flat over $R'$, 
$R'$ is also very strongly $F$-regular by
Lemma~\ref{pure-very-strongly-F-regular.thm}.
\qed

\begin{corollary}
Let $R$ be a noetherian ring of characteristic $p$, and
$\Spec R=\bigcup_i \Spec R_i$ an affine open covering.
Then $R$ is very strongly $F$-regular if and only if $R_i$ is
very strongly $F$-regular for each $i$.
\end{corollary}

\proof Note that each $R_i$ is noetherian.
If $R$ is very strongly $F$-regular, then each $R_i$ is
very strongly $F$-regular by Lemma~\ref{affine-open.thm}.
Conversely, assume that each $R_i$ is very strongly $F$-regular for
each $i$.
Then we can take a finite open subcovering $\Spec R=\bigcup_{j=1}^r
\Spec R_{i_j}$.
As $R'=\prod_{j=1}^r R_{i_j}$ is very strongly $F$-regular by
Lemma~\ref{direct-product.thm} and $R'$ is faithfully flat over 
$R$, $R$ is also very strongly $F$-regular by
Lemma~\ref{pure-very-strongly-F-regular.thm}.
\qed

\begin{lemma}\label{pure-sfr.thm}
Let $R\rightarrow S$ be a cyclically pure ring homomorphism
between noetherian rings of characteristic $p$.
If $S$ is strongly $F$-regular, then so is $R$.
\end{lemma}

\proof Note that $R$ is pure in $S$.
Let $\fm\in\Max R$.
Then there exists some $\fn\in \Max S$ such that $\fn$ lies on $\fm $
and $R_\fm\rightarrow S_\fn$ is pure, see \cite[(9.10)]{HO}.
As $S_\fn$ is very strongly $F$-regular by Lemma~\ref{strongly-F-regular.thm},
$R_\fm$ is very strongly $F$-regular by
Lemma~\ref{pure-very-strongly-F-regular.thm}.
So $R$ is strongly $F$-regular by Lemma~\ref{strongly-F-regular.thm} again.
\qed

\paragraph Let $(R,\fm,K)$ be a complete noetherian 
local ring with a coefficient
field $K\subset R$ of characteristic $p$.
We fix a $p$-base $\Lambda$ of $K$.
A subset $\Gamma$ of $\Lambda$ is said to be cofinite if
$\Lambda\setminus\Gamma$ is a finite set.
For a cofinite subset $\Gamma$ of $\Lambda$ and $e\geq 0$, 
we denote $K^\Gamma_e$ the extension field of $K$ generated by
the all $p^e$th roots of elements in $\Gamma$.
We define $R^\Gamma_e$ to be the completion of the noetherian local ring
$K_e^\Gamma\otimes_K R$.
We denote the inductive limit $\indlim_e R^\Gamma_e $ by $R^\Gamma$.
For an $R$-algebra $A$ essentially of finite type, we define 
$A^\Gamma:=A\otimes_RR^\Gamma$. 

\begin{lemma}\label{ci.thm}
The canonical map $A\rightarrow A^\Gamma$ is a faithfully
flat homomorphism with complete intersection fibers.
\end{lemma}

\proof We may assume that $A$ is a field.
Note that $R\rightarrow K^\Gamma_e\otimes_K R$
is a faithfully flat map with complete intersection fibers
(to verify this, we may assume that $R=K[[x_1,\ldots,x_n]]$).
As $R$ is complete, $K^\Gamma_e\otimes_K R$ is a homomorphic image of a 
regular local ring, and hence the completion
$K^\Gamma_e\otimes_K R\rightarrow R^\Gamma_e$ has complete intersection fibers.

Thus $A^\Gamma=A\otimes_R \indlim R^\Gamma_e\cong\indlim A\otimes_R R^\Gamma_e$
is a noetherian inductive limit of artinian local rings and 
faithfully flat purely inseparable complete intersection homomorphisms.
As $A^\Gamma$ is noetherian, and 
an element of $A^\Gamma$ is either a unit or nilpotent, 
$A^\Gamma$ is artinian local.
The maximal ideal $\fm^\Gamma$ of $A^\Gamma$ is generated by
the maximal ideal $\fm_e^\Gamma$ of the artinian local ring $A^\Gamma_e:=
A\otimes_R R_e^\Gamma$
for sufficiently large $e$.
Then $A^\Gamma_e$ is a complete intersection, and the fiber of
$A^\Gamma_e\rightarrow A^\Gamma$ is a field.
So $A^\Gamma$ is a complete intersection, as desired.
\qed

\begin{lemma}\label{intersection-K.thm}
Set $K^\Gamma=\bigcup_e K^\Gamma_e$.
Then $\bigcap_\Gamma K^\Gamma=K$.
\end{lemma}

\proof Note that $K^\Gamma$ has a basis
\begin{multline*}
B^\Gamma:=\{\xi_1^{\lambda_1}\cdots \xi_n^{\lambda_n}\mid n\geq 0,\;
\xi_1,\ldots,\xi_n 
\text{ are distinct elements in $\Lambda$},
\\
0<\lambda_i<1,\;
\lambda_i\in \Bbb Z[1/p]\}.
\end{multline*}
Thus a linear combination of elements of $B^\Lambda$ lies in
$K^\Gamma$ if and only if any basis element with a nonzero coefficient
lies in $B^\Gamma$.
Thus $\bigcap_\Gamma K^\Gamma$ has a basis $\{1\}$, and hence
$\bigcap_\Gamma K^\Gamma=K$.
\qed

\begin{lemma}\label{intersection-trivial.thm}
Let $K\subset L$ be a field extension, and $\{K_\lambda\}$ a family of
intermediate fields.
Assume that $\bigcap_\lambda K_\lambda=K$.
Let $V$ be a $K$-vector space.
Then $\bigcap_\lambda (K_\lambda\otimes_K V)=V$.
\end{lemma}

\proof
Let $A$ be a basis of $V$ over $K$, and $\sum_{\alpha\in A} c_\alpha \alpha$
be an element of $L\otimes_K V$ with $c_\alpha\in L$.
It lies in $K_\lambda\otimes_K V$ if and only if for any $\alpha$, 
$c_\alpha\in K_\lambda$.
So it lies in $\bigcap_\lambda (K_\lambda\otimes_K V)$ if and only if 
for any $\alpha$, $c_\alpha\in K$.
Thus $\bigcap_\lambda (K_\lambda\otimes_K V)=V$.
\qed

\begin{lemma}\label{R-intersection.thm}
Let $(R,\fm)$ be a noetherian complete local ring of characteristic $p$
with a coefficient field $K$.
Let $\Lambda$ be a $p$-base of $K$.
Then $\bigcap_\Gamma R^\Gamma=R$, where the intersection is taken over
the all cofinite subsets $\Gamma$ of $\Lambda$.
\end{lemma}

\proof First consider the case that $R$ is artinian.
Then $R^\Gamma_e=K^\Gamma_e \otimes_K R$ (completion is unnecessary, since
$K^\Gamma_e \otimes_K R$ is already complete).
So $R^\Gamma=K^\Gamma\otimes_K R$.
So $\bigcap_\Gamma R^\Gamma=\bigcap_\Gamma(K^\Gamma\otimes_K R)=R$
by Lemma~\ref{intersection-K.thm} and Lemma~\ref{intersection-trivial.thm}.

Next  consider the general case.
Let $a$ be an element of $\bigcap_\Gamma R^\Gamma$.
Then $a$ modulo $\fm^n R^\Lambda$ lies in $\bigcap_\Gamma (R/\fm^n)^\Gamma
=R/\fm^n$.
So $a$ lies in $\projlim R/\fm^n=R$.
\qed

\begin{lemma}\label{reg.thm}
Let $R$ be as above, and
$\pi^\Gamma:\Spec A^\Gamma\rightarrow\Spec A$ be the canonical 
homeomorphism.
Then there exists some cofinite subset $\Gamma_0\subset \Lambda$ such
that for every cofinite subset $\Gamma\subset\Gamma_0$, 
$\Reg(A^\Gamma)=(\pi^\Gamma)^{-1}(\Reg(A))$.
Similarly for complete intersection, Gorenstein, Cohen--Macaulay,
$(S_i)$, $(R_i)$, normal, and reduced loci.
\end{lemma}

\proof We prove the assertion for the regular property.
Set $Z^\Gamma$ to be $\pi^\Gamma(\Sing(A^\Gamma))$.
By Kunz's theorem \cite{Kunz}, $Z_\Gamma$ is closed.
The set $\{Z^\Gamma\mid \Gamma:\text{ cofinite in }\Lambda\}$ has 
a minimal element $Z^{\Gamma_0}$.
For any $\Gamma\subset\Gamma_0$, $\Sing R\subset Z^\Gamma= Z^{\Gamma_0}$.
So it suffices to show that $Z^{\Gamma_0}\subset \Sing R$.
Assume the contrary, and let $P\in Z^{\Gamma_0}\setminus\Sing R$.
By \cite[(6.13)]{HH2}, there exists some $\Gamma\subset\Gamma_0$ such
that $PA^\Gamma$ is a prime and $\kappa(P)\otimes_A A^\Gamma$ is a field.
As $A_P$ is regular, $A^\Gamma_{PA^\Gamma}$ is regular.
So $P\notin Z^\Gamma=Z^{\Gamma_0}$.
A contradiction.

The assertion for $(R_i)$ follows easily from this.

The assertions for complete intersection, Gorenstein, Cohen--Macaulay, and
$(S_i)$ properties are true for $\Gamma=\Lambda$, because
$\pi^\Lambda$ is a faithfully flat homomorphism with complete intersection
fibers by Lemma~\ref{ci.thm}.

The normality is $(R_1)+(S_2)$, and reduced property is $(R_0)+(S_1)$.
They are 
proved easily from the assertions for $(R_i)$ and $(S_i)$.
\qed

\begin{lemma}\label{intersection-field.thm}
Let $R$ be as above.
Let $A$ be a field which is essentially of finite type over $R$.
Then $\bigcap_\Gamma A^\Gamma=A$.
\end{lemma}

\proof Let $P$ be the kernel of $R\rightarrow A$.
As $A$ is a field, $P$ is a prime ideal.
Replacing $R$ by $R/P$, 
we may assume that $R$ is a domain and $R\rightarrow A$ is injective.
In view of 
Lemma~\ref{intersection-trivial.thm}, 
replacing $A$ by the field of fractions $Q(R)$ 
of $R$, we may assume that
$A=Q(R)$.

Let $B$ be the normalization of $R$.
Take a cofinite subset $\Gamma_0\subset\Lambda$ such that $R^\Gamma$ is a
domain and $B^{\Gamma}$ is
a normal domain for $\Gamma\subset\Gamma_0$.
This is possible by \cite[(6.13)]{HH2} and Lemma~\ref{reg.thm}.
As $Q(R)\rightarrow Q(R)\otimes_R B$ is an isomorphism, 
$Q(R^\Gamma)=Q(R)\otimes_R R^{\Gamma}=Q(R)\otimes_R B^{\Gamma}$.
Hence $R^\Gamma\rightarrow B^\Gamma$ is the normalization.
In particular, $Q(R^\Gamma)=Q(B^\Gamma)=Q(R)^\Gamma$.

Note that for $\Gamma\subset\Gamma_0$, 
$B^{\Gamma_0}\cap Q(R^\Gamma)$ is purely inseparable, 
hence is integral, over $B^\Gamma$.
Hence $B^{\Gamma_0}\cap Q(B^\Gamma)=B^\Gamma$.

Let $d\in\bigcap_\Gamma B^\Gamma$.
If $c\neq 0$ 
is an element of the conductor $R:B$, then $c\in R^\Gamma:B^\Gamma$.
So $cd\in\bigcap_\Gamma R^\Gamma=R$ by Lemma~\ref{R-intersection.thm}.
So $d\in Q(R)$.
As $d$ is integral over $B$ and $B$ is normal, $d\in B$.
This shows that $\bigcap_\Gamma B^\Gamma=B$.

Now let $\alpha\in\bigcap_\Gamma Q(R)^\Gamma$.
Then there exists some $a\in R^\circ$ and $b\in R^{\Gamma_0}$ such that
$\alpha=b/a$.
As $\alpha\in\bigcap_\Gamma Q(R)^\Gamma$ and $a\in R^\circ$, 
\[
b=a\alpha\in \bigcap_\Gamma(Q(R)^\Gamma\cap R^{\Gamma_0})
\subset\bigcap_\Gamma B^\Gamma=B.
\]
Hence $\alpha=b/a\in Q(R)$.
So $\bigcap_\Gamma Q(R)^\Gamma=Q(R)$.
This is what we wanted to prove.
\qed

\begin{lemma}\label{E.thm}
Let $(R,\fm)$ be a noetherian local ring, and $M$ an $R$-module.
If $\Supp M=\{\fm\}$, $\Ext^1_R(R/\fm,M)=0$, 
and
$\Hom_R(R/\fm,M)\cong R/\fm$,
then $M$ is isomorphic to the injective hull of $R/\fm$.
\end{lemma}

\proof As $\Supp M=\{\fm\}$, $M$ is an essential extension of $\Hom_R(R/\fm
,M)\cong R/\fm$.
So there is an exact sequence of the form
\[
0\rightarrow M\rightarrow E \rightarrow W\rightarrow 0,
\]
where $E$ is the injective hull of $R/\fm$.
As $\Ext^1_R(R/\fm,M)=0$,
we have that
\[
0\rightarrow\Hom_R(R/\fm,M)\rightarrow\Hom_R(R/\fm,E)\rightarrow
\Hom_R(R/\fm,W)\rightarrow 0
\]
is exact.
So $\Hom_R(R/\fm,W)=0$.
As $\Supp W\subset\{\fm\}$, $W=0$.
This shows that $M\cong E$.
\qed

\begin{lemma}\label{E2.thm}
Let $\varphi:(R,\fm)\rightarrow (S,\fn)$ be a local homomorphism
between noetherian local rings of characteristic $p$.
Assume that $\varphi$ is flat, and $S/\fm S$ is Gorenstein of dimension zero.
Then 
\begin{description}
\item[1] $E_R\otimes_R S\cong E_S$, where $E_R$ and $E_S$ respectively denote
the injective hulls of the residue field of $R$ and $S$.
\item[2] Assume further that $\fn=\fm S$.
If $c\in R$, $e\geq 0$, and $c^{(-e)}F^e_R:R\rightarrow R^{(-e)}$ is
$R$-pure, then $c^{(-e)}F^e_S:S\rightarrow S^{(-e)}$ is $S$-pure.
\end{description}
\end{lemma}

\proof {\bf 1}
It is easy to see that $\Supp E_R\otimes_R S=\{\fn\}$.

There is a spectral sequence
\[
E^{p,q}_2=\Ext^p_{S/\fm S}(S/\fn,\Ext^q_S(S/\fm S,E_R\otimes_R S))
\Rightarrow \Ext^{p+q}_S(S/\fn,E_R\otimes_R S).
\]
Note that
\[
\Ext^q_S(S/\fm S,E_R\otimes_R S)\cong
\Ext^q_R(R/\fm,E_R)\otimes_R S\cong
\begin{cases}
S/\fm S & (q=0)\\
0 & (q\neq 0)
\end{cases}.
\]
So $E^{p,q}_2=0$ for $q\neq 0$.
As $S/\fm S$ is the injective hull of the residue field of $S/\fm S$, 
$E^{p,0}_2=0$ for $p>0$, and $E^{0,0}_2\cong S/\fn$.
By Lemma~\ref{E.thm}, 
$E_R\otimes_R S\cong E_S$.

{\bf 2}
Let $\xi$ be a generator of the socle of $E_R$.
Then $\xi\otimes 1\in E_R\otimes_R S$ generates a submodule isomorphic to
$R/\fm\otimes_R S\cong S/\fn$.
Thus $\xi\otimes 1$ is a generator of the socle of $E_S$.
Consider the commutative diagram
\[
\xymatrix{
E_R\otimes_R S~~ \ar[r]^-{1\otimes c^{(-e)}F^e_S} &
~~(E_R\otimes_R S)\otimes_S S^{(-e)}=E_R\otimes_R S^{(-e)} \\
E_R \ar[u]_\varphi \ar[r]^{1\otimes c^{(-e)}F^e_R}
&
E_R\otimes_R R^{(-e)}
\ar[u]_{1\otimes\varphi^{(-e)}}
}.
\]
Then $\xi\in E_R$ goes to a nonzero element in $E_R\otimes_R S^{(-e)}$,
since $c^{(-e)}F^e_R$ is $R$-pure, and $\varphi^{(-e)}$ is faithfully flat.
Thus the socle element $\xi\otimes 1\in E_R\otimes_R S$ goes to
a nonzero element by $1\otimes c^{(-e)}F^e_S$.
This shows that $c^{(-e)}F^e_S: S\rightarrow S^{(-e)}$ is $S$-pure.
\qed

\begin{lemma}\label{excellent.thm}
Let $R$ be an excellent local ring of characteristic $p$,
and $A$ an $R$-algebra essentially of finite type.
Let $c\in A$ such that $A[1/c]$ is regular.
If $c^{(-e)}F^e:A\rightarrow A^{(-e)}$ is $A$-pure for some $e\geq 0$, then
$A$ is very strongly $F$-regular.
\end{lemma}

\proof Let $\hat R$ be the completion of $R$, and $\hat A:=\hat R\otimes_R
A$.
As $R$ is excellent, $\hat A[1/c]$ is regular.
Note that $c^{(-e)}F^e_{\hat A}: \hat A \rightarrow \hat A^{(-e)}$ is
the composite
\[
\hat A 
\xrightarrow{1\otimes c^{(-e)}F^e_A}
\hat A\otimes_A A^{(-e)}
\xrightarrow{\Psi_e(A,\hat A)}
\hat A^{(-e)}.
\]
As $A\rightarrow \hat A$ is regular, this is $\hat A$-pure.
Replacing $R$ by its completion $\hat R$ and $A$ by $\hat A$,
we may assume that $R$ is complete local by
Lemma~\ref{pure-very-strongly-F-regular.thm}.

Take a coefficient field $K$ of $R$, fix a $p$-base $\Lambda$ of $K$, and
take a cofinite subset $\Gamma_0$ of $\Lambda$ such that $A^\Gamma[1/c]$ is
regular for any cofinite subset $\Gamma\subset\Gamma_0$.

For a cofinite subset $\Gamma$ of $\Lambda$, let $\pi^\Gamma:\Spec A^\Gamma
\rightarrow \Spec A$ be the canonical map.
Let $W^\Gamma$ be the closed subset of $\Spec A^\Gamma$ consisting of
prime ideals $P$ such that $c^{(-e)}F_{A_P^\Gamma}^e:A_P^\Gamma
\rightarrow
(A_P^\Gamma)^{(-e)}$ is not $A_P^\Gamma$-pure
(it is closed, since $A^\Gamma$ is $F$-finite by \cite[(6.6), (6.8)]{HH2}).
Let $Z^\Gamma=\pi^\Gamma(W^\Gamma)$.
It is easy to see that if $\Gamma'\subset\Gamma$, then $Z^{\Gamma'}
\subset Z^{\Gamma}$.
Let $\Gamma_1\subset\Gamma_0$ be a cofinite subset such that
$Z^{\Gamma_1}$ is minimal.
We show that $Z^{\Gamma_1}$ is empty.
Assume the contrary.
Then there is a prime ideal $P\in Z^{\Gamma_1}$.
Take $\Gamma_2\subset \Gamma_1$ such that $PA^{\Gamma_2}$ is a prime ideal.
This can be done by \cite[(6.13)]{HH2}.
As $c^{(-e)}F^e_{A_P}:A_P\rightarrow A_P^{(-e)}$ is $A_P$ pure,
$c^{(-e)}F^e_{A^{\Gamma_2}_P}:A_{P}^{\Gamma_2}\rightarrow 
(A_P^{\Gamma_2})^{(-e)}$ is $A_P^{\Gamma_2}$-pure by Lemma~\ref{E2.thm}.
On the other hand, as 
$P\in Z^{\Gamma_1}=Z^{\Gamma_2}$, $PA^{\Gamma_2}\in W^{\Gamma_2}$.
A contradiction.
So $c^{(-e)}F^e_{A^{\Gamma_1}}:A^{\Gamma_1}\rightarrow (A^{\Gamma_1})^{(-e)}$ 
is $A^{\Gamma_1}$-pure.

As $A^{\Gamma_1}$ is $F$-finite \cite[(6.6), (6.8)]{HH2}, 
$A^{\Gamma_1}[1/c]$ is
regular, and $c^{(-e)}F^e:A^{\Gamma_1}\rightarrow (A^{\Gamma_1})^{(-e)}$ is 
$A^{\Gamma_1}$-pure
(in particular, $c$ is a nonzerodivisor of $A^{\Gamma_1}$),
$A^{\Gamma_1}$ is very strongly $F$-regular by \cite[Theorem~3.3]{HH3}.
Since $A\rightarrow A^{\Gamma_1}$ is faithfully flat, $A$ is very strongly
$F$-regular.
\qed

\begin{lemma}\label{smooth-base-change.thm}
Let $A\rightarrow B$ be a regular homomorphism of noetherian rings of
characteristic $p$.
Assume that $A$ is very strongly $F$-regular \(resp.\ strongly $F$-regular\), 
and is excellent.
Assume also that $B$ is essentially of finite type over an excellent
local ring \(resp.\ locally excellent\).
Then $B$ is very strongly $F$-regular \(resp.\ strongly $F$-regular\).
\end{lemma}

\proof First consider the case that $A$ is very strongly $F$-regular.
Take $c\in A^\circ$ such that $A[1/c]$ is regular.
This is possible, since $A$ is normal and excellent.
Then we can take $e>0$ such that $c^{(-e)}F^e:A\rightarrow A^{(-e)}$ is
$A$-pure.
Then plainly, $1_B\otimes c^{(-e)}F^e:B\rightarrow B\otimes_A A^{(-e)}$
is $B$-pure.
As $\Psi_e:B\otimes_A A^{(-e)}\rightarrow B^{(-e)}$ is faithfully flat
by Radu \cite{Radu} and Andr\'e \cite{Andre}, 
$c^{(-e)}F^e_B:B\rightarrow B^{(-e)}$ is $B$-pure.
As $B[1/c]$ is regular, $B$ is very strongly $F$-regular by
Lemma~\ref{excellent.thm}.

Next consider the case that $A$ is strongly $F$-regular.
Take $\fn\in\Max B$, and set $\fm:=A\cap \fn$.
Then $A_\fm$ is very strongly $F$-regular by
Lemma~\ref{strongly-F-regular.thm} and Lemma~\ref{vsfr-localization.thm}.
By the first paragraph, $B_\fn$ is very strongly $F$-regular.
By Lemma~\ref{strongly-F-regular.thm}, $B$ is strongly $F$-regular.
\qed

For a noetherian ring $A$ of characteristic $p$, 
set
\begin{eqnarray*}
\SFR(A) &:= &\{P\in\Spec A\mid A_P\text{ is strongly $F$-regular}\},\\
\NonSFR(A) &:= &\Spec A\setminus\SFR(A).
\end{eqnarray*}

\begin{lemma}\label{open-F-finite.thm}
Let $A$ be an $F$-finite noetherian ring of characteristic $p$.
Then $\SFR(A)$ is a Zariski open subset of $\Spec A$.
\end{lemma}

\proof As the reduced locus of $A$ is open, 
we may assume that $A$ is reduced.
Take $c\in A^\circ$ such that $A[1/c]$ is regular.
For each $e\geq 0$, 
let $U_e$ be the complement of the support of the cokernel of the map
\[
c^{(-e)}F^e:\Hom_A(A^{(-e)},A)\rightarrow \Hom_A(A,A)=A
\qquad
(\varphi\mapsto \varphi c^{(-e)}F^e).
\]
Then $U_e$ is open, and $\SFR(A)=\bigcup_e U_e$ is also open.
\qed

\begin{lemma}\label{zariski-open.thm}
Let $(R,\fm)$ be a complete local ring with the coefficient field $K$,
and $A$ an $R$-algebra essentially of finite type.
Let $\Lambda$ be a $p$-base of $K$.
For a cofinite subset $\Gamma\subset\Lambda$, 
let $\pi^\Gamma:\Spec A^\Gamma\rightarrow\Spec A$ be the canonical map.
Then there exists some cofinite subset $\Gamma_0$ of $\Lambda$ such that
$\pi^\Gamma(\SFR(A^\Gamma))=\SFR(A)$ for any cofinite subset $\Gamma$ of
$\Gamma_0$.
In particular, $\SFR(A)$ is a Zariski open subset of $\Spec A$.
\end{lemma}

\proof 
Let $Z^\Gamma$ be the closed subset $\pi^\Gamma(\NonSFR(A^\Gamma))$
of $\Spec A$.
Take $\Gamma_0$ such that $Z^{\Gamma_0}$ is minimal.
Then $Z^\Gamma=Z^{\Gamma_0}\supset \NonSFR(A)$ for any cofinite subset
$\Gamma\subset\Gamma_0$.
Assume that $P\in Z^{\Gamma_0}\setminus \NonSFR(A)$.
By \cite[(6.13)]{HH2}, we can take
$\Gamma_1\subset\Gamma_0$ such that $PA^{\Gamma}$ is a prime ideal
for any cofinite subset $\Gamma\subset\Gamma_1$.
Take $c\in A_P^\circ$ such that $A_P[1/c]$ is regular.
We can take $\Gamma\subset \Gamma_1$ such that $A^\Gamma_P[1/c]\cong 
A_P[1/c]\otimes_R R^\Gamma$ is regular by Lemma~\ref{reg.thm}.
As $A_P$ is very strongly $F$-regular,
there exists some $e>0$ such that $c^{(-e)}F^e_{A_P}:A_P\rightarrow
A_P^{(-e)}$ is $A_P$-pure.
By Lemma~\ref{E2.thm}, $c^{(-e)}F^e_{A_P^\Gamma}:A_P^\Gamma
\rightarrow (A_P^\Gamma)^{(-e)}$ is $A_P^\Gamma$-pure.
As $A_P^\Gamma[1/c]$ is regular, 
$A_P^\Gamma$ is strongly $F$-regular.
This contradicts the choice of $P$.
So $\pi^\Gamma(\SFR(A^\Gamma))=\SFR(A)$ for $\Gamma\subset\Gamma_0$,
as desired.

Now the openness of $\SFR(A)$ follows from
Lemma~\ref{open-F-finite.thm}.
\qed

\begin{corollary}\label{Gamma-sfr.thm}
Let $A$ be as in {\rm Lemma~\ref{zariski-open.thm}}.
Assume that $A$ is strongly $F$-regular.
Then there exists some cofinite subset $\Gamma_0$ of $\Lambda$ such that
for any cofinite subset $\Gamma$ of $\Gamma_0$, 
$A^\Gamma$ is strongly $F$-regular.
\qed
\end{corollary}

\begin{lemma}\label{closed.thm}
Let $\varphi:X\rightarrow Y$ be a continuous map between topological
spaces, and $Z\subset X$ a closed subset.
Assume that $X$ is a noetherian topological space, 
and each irreducible closed subset of $X$ has a generic point.
If $\varphi(Z)$ is
closed under specialization,
then $\varphi(Z)\subset Y$ is closed.
\end{lemma}

\proof By assumption, there is a finite set of points 
$C=\{z_1,\ldots,z_r\}$ of $Z$ such that
the closure of $C$ is $Z$.
As $\varphi(C)\subset\varphi(Z)$ and $\varphi(Z)$ is closed under
specialization, we have $\overline{\varphi(C)}=\bigcup_i \overline{
\{\varphi(z_i)\}}\subset\varphi(Z)$.
As $\varphi$ is continuous, 
\[
\varphi(Z)=\varphi(\bar C)\subset\overline{\varphi(C)}\subset\varphi(Z).
\]
Hence $\varphi(Z)=\overline{\varphi(C)}$ is closed.
\qed

\begin{proposition}\label{open.thm}
Let $R$ be an excellent local ring of characteristic $p$, and
$A$ an $R$-algebra essentially of finite type.
Then $\SFR(A)$ is Zariski open in $\Spec A$.
\end{proposition}

Let $\hat R$ be the completion of $R$, and set $\hat A:=\hat R\otimes_R A$.
Let $\rho:\Spec \hat A\rightarrow \Spec A$ be the map associated with the
base change of the completion.
Let $Q$ be a prime ideal of $\hat A$.
Then $\hat A_Q$ is strongly $F$-regular if and only if so is $A_{\rho(Q)}$
by Lemma~\ref{pure-sfr.thm} and Lemma~\ref{smooth-base-change.thm}.
Thus $\rho^{-1}(\NonSFR(A))=\NonSFR(\hat A)$.
Letting $X=\Spec \hat A$, $Y=\Spec A$, and $Z=\NonSFR(\hat A)$, 
$\rho(Z)=\rho(\rho^{-1}(\NonSFR(A)))=\NonSFR(A)$ is closed
by Lemma~\ref{closed.thm}, since it is closed under specialization
by Lemma~\ref{strongly-F-regular.thm}.
\qed

\begin{proposition}\label{excellent2.thm}
Let $R$ be an excellent noetherian local ring of characteristic $p$,
and $A$ an $R$-algebra essentially of finite type.
Let $c\in A$ such that $A[1/c]$ is strongly $F$-regular.
If $c^{(-e)}F^e_A:A\rightarrow A^{(-e)}$ is $A$-pure for some $e\geq 0$,
then $A$ is very strongly $F$-regular.
\end{proposition}

\proof Let $\hat R$ be the completion of $R$, and set $\hat A:=\hat R
\otimes_R A$.
Then $\hat A[1/c]$ is strongly $F$-regular by 
Lemma~\ref{smooth-base-change.thm}.
Moreover, $c^{(-e)}F^e_{\hat A}: \hat A\rightarrow \hat A^{(-e)}$ is 
$\hat A$-pure as in the proof of Lemma~\ref{excellent.thm}.
By Lemma~\ref{pure-very-strongly-F-regular.thm}, 
we may assume that $R$ is complete local.

Now take a coefficient field $K$ of $R$, and take a $p$-base $\Lambda$ of
$K$.
Then by Corollary~\ref{Gamma-sfr.thm}, there exists some cofinite
subset $\Gamma_0$ of $\Lambda$ such that for each cofinite subset
$\Gamma$ of $\Gamma_0$, $A^\Gamma[1/c]=A[1/c]^\Gamma$ is strongly $F$-regular.
As in the proof of Lemma~\ref{excellent.thm},
there exists some $\Gamma_1\subset\Gamma_0$ such that
$c^{(-e)}F^e_{A^{\Gamma_1}}:A^{\Gamma_1}\rightarrow(A^{\Gamma_1})^{(-e)}$ 
is $A^{\Gamma_1}$-pure.
As $A^{\Gamma_1}$ is $F$-finite, $A^{\Gamma_1}$ is very 
strongly $F$-regular by \cite[(3.3)]{HH3}.
By Lemma~\ref{pure-very-strongly-F-regular.thm}, 
$A$ is very strongly $F$-regular.
\qed

\begin{corollary}
Let $R$ be an excellent noetherian local ring of characteristic $p$,
and $A$ an $R$-algebra essentially of finite type.
Then $A$ is very strongly $F$-regular if and only if it is strongly
$F$-regular.
\end{corollary}

\proof By Lemma~\ref{vsfr-localization.thm}, a very strongly $F$-regular
implies strongly $F$-regular.

Conversely, assume that $A$ is strongly $F$-regular.
Then letting $c=1$, $A[1/c]$ is strongly $F$-regular, and
$c^{(0)}F^0_A:A\rightarrow A^{(0)}$ is $A$-pure,
since it is the identity map.
By Lemma~\ref{excellent2.thm}, $A$ is very strongly $F$-regular.
\qed

\begin{corollary}\label{sfr-check.thm}
Let $A$ be an locally excellent noetherian ring of characteristic $p$.
Let $c\in A$ such that $A[1/c]$ is strongly $F$-regular.
If $c^{(-e)}F^e_A:A\rightarrow A^{(-e)}$ is $A$-pure for some $e\geq 0$,
then $A$ is strongly $F$-regular.
\end{corollary}

\proof We may assume that $A$ is local.
Then the assertion is obvious by Lemma~\ref{excellent2.thm}.
\qed

\begin{theorem}\label{Fedder-Watanabe.thm}
Let $\varphi:A\rightarrow B$ be a homomorphism of noetherian rings of
characteristic $p$.
Assume that $A$ is a strongly $F$-regular domain.
Assume that the generic fiber $Q(A)\otimes_A B$ is strongly $F$-regular,
where $Q(A)$ is the field  of fractions of $A$.
If $\varphi$ is $F$-pure and $B$ is locally excellent, then
$B$ is strongly $F$-regular.
\end{theorem}

\proof Assume the contrary.
Then there is a prime ideal $P$ of $B$ such that $B_P$ is not
strongly $F$-regular, but $B_Q$ is strongly $F$-regular for 
any prime ideal
$Q\subsetneq P$.
Replacing $B$ by $B_P$ and $A$ by $A_{P\cap A}$,
we may assume that $(A,\fm)$ and $(B,\fn)$ are local and
$\varphi$ is local, and we may assume that $\NonSFR(B)=\{\fn\}$.
By assumption, $A\neq Q(A)$.
So there is a nonzero element $c\in\fm$.
Then $B[1/c]$ is strongly $F$-regular, since $c\in\fn$.
As $A$ is a very strongly $F$-regular domain, 
there exists some $e\geq 0$ such that
$c^{(-e)}F^e_A:A\rightarrow A^{(-e)}$ is $A$-pure.
As $c^{(-e)}F^e_B:B\rightarrow B^{(-e)}$ is the composite
\[
B=B\otimes_AA\xrightarrow{1_B\otimes c^{(-e)}F^e_A}
B\otimes_A A^{(-e)}
\xrightarrow{\Psi_e(A,B)}
B^{(-e)}
\]
and $\Psi_e(A,B)$ is pure, $c^{(-e)}F^e_B$ is $B$-pure.
By Corollary~\ref{sfr-check.thm}, $B$ is strongly $F$-regular.
This is a contradiction.
\qed

\section{Cohen--Macaulay $F$-injective property}

\paragraph We say that a noetherian local ring $(R,\fm)$ of characteristic $p$
is Cohen--Macaulay $F$-injective (CMFI for short) 
if $R$ is Cohen--Macaulay, and the
Frobenius map $H^d_\fm(R)\rightarrow H^d_\fm(R)\otimes_R R^{(-e)}
\cong H^d_{\fm^{(-e)}}(R^{(-e)})$ is
injective for some (or equivalently, any) $e>0$, where $d$ is the
dimension of $R$.
Obviously, $R$ is CMFI if and only if its completion is.
A noetherian ring of characteristic $p$ is said to be CMFI if
its localization at any maximal ideal is CMFI.

\paragraph
Let $I$ be an ideal of a noetherian ring $R$ of characteristic $p$.
The Frobenius closure of $I$ is defined to be
\[
I^F:=\{x\in R\mid x^q\in I^{[q]}\text{ for some $q=p^e$}\},
\]
where $I^{[q]}=I^{(e)}R$.
If $I=I^F$, then we say that $I$ is Frobenius closed.

\paragraph Let $(R,\fm)$ be a Cohen--Macaulay local ring of characteristic $p$
with dimension $d$.
Take a system of parameters $x_1,\ldots,x_d$ of $R$.
We see that the local cohomology $H^d_\fm(R)$ is the $d$th cohomology
group of the modified \v Cech complex, see \cite[(3.5)]{BH}.
It is identified with the inductive limit of the inductive system:
\begin{equation}\label{lc.eq}
R/(x_1,\ldots,x_d)\xrightarrow{x_1\cdots x_d} R/(x_1^2,\ldots,x_d^2)
\xrightarrow{x_1\cdots x_d} R/(x_1^3,\ldots,x_d^3)\rightarrow\cdots,
\end{equation}
where $a\in R/(x_1^t,\ldots,x_d^t)$ corresponds to $a/(x_1\cdots x_d)^t$ in
$H^d_\fm(R)$.
Note that the maps of the inductive system (\ref{lc.eq}) are injective, because
$x_1,\ldots,x_d$ is a regular sequence.
In particular, $R/(x_1^t,\ldots,x_d^t)$ can be identified with a
submodule of $H^d_\fm(R)$.

So if $R$ is CMFI and $F:H^d_\fm(R)\rightarrow H^d_\fm(R)\otimes_R R^{(-e)}$ is
injective, then $F:R/(x_1,\ldots,x_d)\rightarrow R/(x_1,\ldots,x_d)\otimes_R
R^{(-e)}$ is injective.
In other words, if $x\in R$ and $x^q\in (x_1^q,\ldots,x_d^q)$, then $x\in
(x_1,\ldots,x_d)$.

On the other hand, in (\ref{lc.eq}), the socle of 
$R/(x_1^t,\ldots,x_d^t)$ is mapped bijectively onto the socle of
$R/(x_1^{t+1},\ldots,x_d^{t+1})$, since the map is injective, and
the dimensions of the socles are equal (they agree with the 
Cohen--Macaulay type of $R$).
Hence $\Soc R/(x_1,\ldots,x_d)\rightarrow \Soc H^d_\fm(R)$ is isomorphic.
This shows that if $F:R/(x_1,\ldots,x_d)\rightarrow R/(x_1,\ldots,x_d)\otimes_R
R^{(-1)}$ is injective, then $R$ is CMFI.
Hence we have

\begin{lemma}\label{CMFI-equiv.thm}
Let $(R,\fm)$ be a Cohen--Macaulay local ring of characteristic $p$.
Then the following are equivalent.
\begin{description}
\item[1] $R$ is CMFI.
\item[2] Any parameter ideal
of $R$ is Frobenius closed, where a parameter ideal means an ideal generated
by a system of parameters.
\item[3] For one system of parameters $x_1,\ldots,x_d$ of $R$, 
$x\in R$, $x^p\in (x_1^p,\ldots,x_d^p)$ implies $x\in (x_1,\ldots,x_d)$.
\end{description}
\qed
\end{lemma}

This lemma is a simplified variation of \cite[Proposition~2.2]{FW}.

\paragraph
If $R$ is CMFI and $(x_1,\ldots,x_d)$ is a system of parameters of $R$,
and $0\leq s\leq d$, then
\begin{multline*}
(x_1,\ldots,x_s)^F=
(\bigcap_{i_{s+1},\ldots,i_d\geq 0}(x_1,\ldots,x_s,x_{s+1}^{i_{s+1}},\ldots,
x_d^{i_d}))^F
\subset\\
\bigcap_{i_{s+1},\ldots,i_d\geq 0}(x_1,\ldots,x_s,x_{s+1}^{i_{s+1}},\ldots,
x_d^{i_d})^F
=\\
\bigcap_{i_{s+1},\ldots,i_d\geq 0}(x_1,\ldots,x_s,x_{s+1}^{i_{s+1}},\ldots,
x_d^{i_d})
=(x_1,\ldots,x_s).
\end{multline*}
Thus, an ideal generated by a regular sequence of $R$ is Frobenius closed.

\begin{lemma}\label{ff-cmfi.thm}
Let $f:(R,\fm)\rightarrow (S,\fn)$ be a flat
local homomorphism between
noetherian local rings of characteristic $p$.
If $S$ is CMFI, then so is $R$.
\end{lemma}

\proof As $S$ is Cohen--Macaulay and $f$ is flat local, $R$ is 
Cohen--Macaulay.
Let $x_1,\ldots,x_d$ be a system of parameters of $R$, and set 
$I:=(x_1,\ldots,x_d)$.
Then $I^F\subset (IS)^F\cap R\subset IS\cap R=I$,
because $IS$ is generated by a regular sequence, and $S$ is CMFI.
\qed

Let $R$ be a noetherian ring of characteristic $p$.
We define
\[
\CMFI(R)=\{P\in\Spec R\mid R_P\text{ is CMFI}\}
\]
and
$\NonCMFI(R):=\Spec R\setminus \CMFI(R)$.

\begin{lemma}\label{mcm.thm}
Let $R$ be a Cohen--Macaulay ring, and $M$ a finite $R$-module.
Then the locus
\[
\MCM(M):=\{P\in\Spec R\mid M_P \text{ is a maximal Cohen--Macaulay 
$R_P$-module}\}
\]
is a Zariski open subset of $\Spec R$.
\end{lemma}

\proof Set $S=R\oplus M$ to be the idealization of $M$.
$S$ is a noetherian ring with the product $(r,m)(r',m')=(rr',rm'+mr')$.
Note that for $P\in\Spec R$, $M_P$ is maximal Cohen--Macaulay (MCM for short) 
if and only if $S_P$ is
an MCM $R$-module if and only if $S_P$ is a Cohen--Macaulay ring.
Note that $\pi:S\rightarrow S/M=R$ induces a homeomorphism
$\pi^{-1}:\Spec R\rightarrow \Spec S$, and $S_P=S_{\pi^{-1}(P)}$.
So it suffices  to show that the Cohen--Macaulay locus of the ring $S$
is open.
This is well-known, as $S$ is finitely generated over a Cohen--Macaulay
ring, see \cite[Exercise~24.2]{CRT}.
\qed

\begin{lemma}\label{F-finite-open.thm}
Let $R$ be an $F$-finite noetherian ring of characteristic $p$.
Then $\CMFI(R)$ is Zariski open in $\Spec R$.
\end{lemma}

\proof Note that $R$ is excellent \cite{Kunz}.
So the Cohen--Macaulay locus of $R$ is Zariski open by Nagata's criterion
\cite{CRT}.
Hence we may assume that $R$ is Cohen--Macaulay.
Let $M$ be the cokernel of the Frobenius map $F_R:R\rightarrow R^{(-1)}$.
Then $R_P$ is CMFI if and only if $M_P$ is a maximal Cohen--Macaulay
$R_P$-module \cite{Hashimoto}.
So $\CMFI(R)=\MCM(M)$ is Zariski open by Lemma~\ref{mcm.thm}.
\qed

\begin{lemma}\label{Gamma-F-injective.thm}
Let $(R,\fm)$ be a complete noetherian local ring of characteristic $p$ with
a coefficient field $K$.
Let $\Lambda$ be a $p$-base of $K$.
Let $A$ be an $R$-algebra essentially of finite type.
Then there exists some cofinite subset $\Gamma_0$ of $\Lambda$ such that
for any cofinite subset $\Gamma$ of $\Gamma_0$, $\pi^\Gamma(\CMFI(A^\Gamma))
=\CMFI(A)$, where $\pi^\Gamma:\Spec A^\Gamma\rightarrow \Spec A$ is the
canonical morphism.
\end{lemma}

\proof Set $Z^\Gamma:=\pi^\Gamma(\NonCMFI(A^\Gamma))$.
The set $\{Z^\Gamma \mid \Gamma 
\text{ is cofinite in } \Lambda\}$ is a non-empty set of 
closed subsets of the noetherian space $\Spec A$.
Take $\Gamma_0$ so that $Z^{\Gamma_0}$ is
minimal.
By Lemma~\ref{ff-cmfi.thm}, it is easy to see that $Z^\Gamma=Z^{\Gamma_0}$
for $\Gamma\subset \Gamma_0$.

Clearly, $\CMFI(A)\supset \pi^\Gamma(\CMFI(A^{\Gamma_0}))
=\Spec A\setminus Z^{\Gamma_0}$ 
by Lemma~\ref{ff-cmfi.thm}.
Hence $Z^{\Gamma_0}\supset \NonCMFI(A)$.
So it suffices to show that
$Z^{\Gamma_0}\subset \NonCMFI(A)$.

Assume the contrary, and take
$P\in Z^{\Gamma_0}\cap \CMFI(A)$.
We can take $\Gamma_1\subset\Gamma_0$ such that $PA^{\Gamma_1}$ 
is a prime ideal
of $A^{\Gamma_1}$.
Set $d:=\height P$, and take a parameter ideal $I=(x_1,\ldots,x_d)$ 
of $A_P$.

Set $k:=\kappa(P)$.
For $\Sigma\subset\Gamma_1$, 
\begin{multline*}
\Soc(A_P^\Sigma/IA_P^\Sigma)
=\Hom_{A_P^\Sigma}(k\otimes_{A_P}A_P^\Sigma,
A_P/I\otimes_{A_P}A_P^\Sigma)
\cong\\
\Hom_{A_P}(k,A_P/I)\otimes_{k}k^\Sigma
=\Soc(A_P/I)\otimes_{k}k^\Sigma,
\end{multline*}
where $k^\Sigma=k\otimes_{A_P}A_P^\Sigma$.
Thus $V_k:=\Soc(A_P/I)$ gives a $k$-structure of
the finite dimensional
$k^{\Gamma_1}$-vector space 
$V:=\Soc(A_P^{\Gamma_1}/IA_P^{\Gamma_1})$.
For $\Sigma\subset \Gamma_1$, we set $V^{\Sigma}:=V_k\otimes_{k}
k^{\Sigma}\cong\Soc(A_P^\Sigma/IA_P^\Sigma)$.

Now consider $M:=\Ker(A_P^{\Gamma_1}/IA_P^{\Gamma_1}
\xrightarrow{F}(A_P^{\Gamma_1})^{(-1)}/(IA_P^{\Gamma_1})^{(-1)})$ and
$E:=M\cap V=\Soc M$.
For each $\Sigma\subset\Gamma_1$, we set $E^{\Sigma}=E\cap V^{\Sigma}$.
For $\Sigma'\subset\Sigma\subset \Gamma_1$, 
the canonical map
$k^{\Sigma}\otimes_{k^{\Sigma'}}E^{\Sigma'}
\rightarrow E^{\Sigma}$ is injective,
and hence $\dim_{k^{\Sigma'}}E^{\Sigma'}\leq 
\dim_{k^{\Sigma}}E^{\Sigma}$.
Take a cofinite subset $\Omega\subset\Gamma_1$ such that 
$\dim_{k^{\Omega}}E^{\Omega}$ is small as possible.
Let $\kappa$ be the smallest field of definition of $E^{\Omega}$ over
$k$.
Namely, $\kappa$ is the smallest intermediate field $k\subset
\kappa\subset k^\Omega$ such that
$k^\Omega\otimes_\kappa (E^\Omega\cap(V_k\otimes_{k}\kappa))
\rightarrow E^\Omega$ is surjective, see \cite[(3.10)]{Nagata}.
Then by the choice of $\Omega$, for any cofinite subset $\Omega'$ of 
$\Omega$, $k^{\Omega'}\supset \kappa$.
Hence $\kappa\subset\bigcap_{\Omega'}k^{\Omega'}=k$ by
Lemma~\ref{intersection-field.thm}.
Hence $\kappa=k$.

As the diagram
\[
\xymatrix{
A_P/I \otimes_{A_P} A^\Sigma_P
\ar[r]^-F &
A_P/I \otimes_{A_P} (A^{\Sigma}_P)^{(-1)} \\
A_P/I
\ar@{_{(}->}[u] \ar[r]^-F &
A_P/I \otimes_{A_P} A_P^{(-1)}
\ar@{_{(}->}[u]
}
\]
is commutative and the bottom $F$ is injective by the assumption
$P\in\CMFI(A)$, $M\cap A_P/I =0$.
In particular, $E_k:=E^\Omega \cap V_k=0$.
As $k$ is the field of definition of $E^\Omega$, $E^\Omega=0$.
This shows that 
$F:A_P^\Omega/IA_P^\Omega\rightarrow (A_P^\Omega)^{-1}/I(A_P^\Omega)^{(-1)}$ 
is injective.
As $A_P^\Omega$ is Cohen--Macaulay, $A_P^\Omega$ is CMFI by
Lemma~\ref{CMFI-equiv.thm}.
This contradicts $P\in Z^{\Gamma_0}=Z^{\Omega}$.
\qed

\begin{corollary}
Let $(R,\fm)$ be a noetherian local ring of characteristic $p$,
and $A$ a finite $R$-algebra.
If $A$ is CMFI, then there exists some faithfully flat $F$-finite 
local $R$-algebra
$R'$ such that $A'=R'\otimes_R A$ is CMFI.
\end{corollary}

\proof Let $\hat R$ be the  completion of $R$.
Then $A$ is a semilocal ring, and $\hat R\otimes_R A$
is the direct product of the completions of the local rings of $A$
at the maximal ideals.
Hence $\hat R\otimes_R A$ is CMFI.
So replacing $R$ by $\hat R$ and $A$ by $\hat R\otimes_R A$, we may
assume that $R$ is complete.

Let $K$ be a coefficient field of $R$, and take a $p$-base $\Lambda$ of
$K$.
Then by Lemma~\ref{Gamma-F-injective.thm}, there exists some cofinite subset
$\Gamma$ such that $A'=A^\Gamma$ is CMFI.
\qed

\begin{corollary}\label{CMFI-localization.thm}
Let $(R,\fm)$ be a noetherian local CMFI ring of characteristic $p$.
Then for any prime ideal $P$ of $R$, $R_P$ is CMFI.
\end{corollary}

\proof Let $(R',\fm')$ 
be an $F$-finite CMFI local $R$-algebra which is faithfully
flat over $R$.
There exists some $Q\in\Spec R'$ such that $Q\cap R=P$.
By Lemma~\ref{F-finite-open.thm}, 
the CMFI locus of $R'$ is open.
As $\fm'\in\CMFI(R')$, 
we have that $\CMFI(R')=\Spec R'$.
Thus $R'_Q$ is CMFI.
By Lemma~\ref{ff-cmfi.thm}, $R_P$ is also CMFI.
\qed

\begin{lemma}\label{artinian-free.thm}
A flat module over an artinian local ring is free.
\end{lemma}

\proof Let $(R,\fm)$ be an artinian local ring, and $F$ a flat $R$-module.
Then by \cite[(III.2.1.8)]{Hashimoto2}, there is a short exact sequence
\[
0\rightarrow P\rightarrow F\rightarrow G\rightarrow0
\]
of $R$-modules in which $P$ is $R$-free, $G$ is $R$-flat, and
$G/\fm G=0$.
Then by \cite[(I.2.1.6)]{Hashimoto2}, $G=0$.
\qed

\begin{corollary}\label{cok-free.thm}
Let $(R,\fm)$ be an artinian local ring,
and $f:P\rightarrow F$ an $R$-linear map between flat $R$-modules.
If $\bar f: P/\fm P\rightarrow F/\fm F$ is injective, then $f$ is $R$-pure
and $\Coker f$ is $R$-free.
\end{corollary}

\proof Follows immediately by Lemma~\ref{artinian-free.thm} 
and \cite[(I.2.1.4)]{Hashimoto2}.
\qed

\begin{corollary}\label{weak-regular-sequence.thm}
Let $(R,\fm)$ be an artinian local ring, $A$ an $R$-algebra, and
$M$ an $R$-flat $A$-module.
Let $(x_1,\ldots,x_n)$ be a sequence in $A$.
If $(x_1,\ldots,x_n)$ is a \(weak\) $M/\fm M$-sequence,
then $(x_1,\ldots,x_n)$ is a \(weak\) 
$M$-sequence, and $M/(x_1,\ldots,x_n)M$ is
$R$-flat.
\end{corollary}

\proof Set $M_0=M$, $M_i=M/(x_1,\ldots,x_i)M$, $\bar M=M/\fm M$, and $\bar M_i
=\bar M/(x_1,\ldots,x_i)\bar M$.
We prove that $M_i$ is $R$-flat and $(x_1,\ldots,x_i)$ is a weak 
$M$-sequence by
induction on $i$.
If $i=0$, then $M$ is $R$-flat by assumption, and the empty sequence is
an $M$-sequence of length zero.
If $i>0$, then $M_{i-1}$ is $R$-flat and $(x_1,\ldots,x_{i-1})$ is an
$M$-sequence by induction assumption.
As $x_i: \bar M_{i-1}\rightarrow \bar M_{i-1}$ is an injective $R$-linear map,
$x_i:M_{i-1}\rightarrow M_{i-1}$ is injective, and $M_i=M_{i-1}/x_i M_{i-1}$
is $R$-flat.
Clearly, if $\bar M_n\neq 0$, then $M_n\neq 0$, and $(x_1,\ldots,x_n)$ is
a regular sequence.
\qed

\begin{lemma}\label{key.thm}
Let $(R,\fm)$ be an artinian local ring, $A$ an $R$-algebra, and
$f:M\rightarrow N$ an $A$-linear map between $R$-flat $A$-modules.
Let $x_1,\ldots,x_n$ be a sequence of elements in $A$, and assume that
$x_1,\ldots,x_n$ is both a weak $M/\fm M$-sequence and a weak 
$N/\fm N$-sequence.
Assume that $\bar f_n:\bar M_n\rightarrow \bar N_n$ is injective,
where
$M_i:=M/(x_1,\ldots,x_i)M$, $N_i:=N/(x_1,\ldots,x_i)N$, 
$\bar M_i=R/\fm\otimes_R M_i$, and
$\bar N_i=R/\fm\otimes_R N_i$.
Then $f_n:M_n\rightarrow N_n$ is injective, and $\Coker f_n$ is $R$-flat.
\end{lemma}

\proof By Corollary~\ref{weak-regular-sequence.thm},
$M_n$ and $N_n$ are $R$-flat.
Since $\bar f_n$ is injective, the assertions follow immediately
by Corollary~\ref{cok-free.thm}.
\qed

The following was first proved by Aberbach--Enescu \cite{AE}, 
see \cite{Enescu2}.

\begin{proposition}\label{base-change-CMFI.thm}
Let $\varphi:(A,\fm)\rightarrow (B,\fn)$ be 
a flat local homomorphism of noetherian local rings of characteristic $p$.
If $A$ is CMFI and $B/\fm B$ is geometrically CMFI over $A/\fm$ 
\(see {\rm\cite[Definition~5.3]{Hashimoto}}\), then $B$ is CMFI.
\end{proposition}

\proof Clearly, $B$ is Cohen--Macaulay.
Take a system of parameters $x_1,\ldots,x_n$ of $A$ and a sequence
$y_1,\ldots,y_m$ in $\fn$ whose image in $B/\fm B$ is a system
of parameters.
Set $I=(x_1,\ldots,x_n)A$, $J=(y_1,\ldots,y_m)B$, and $\fa=IB+J$.
$B/J$ is $A$-flat by \cite[Corollary to (22.5)]{CRT}.
Since $F_A:A/I\rightarrow A/I\otimes_A A^{(-1)}$ is injective by the
CMFI property of $A$, 
$F_A: B/J\otimes_A A/I\rightarrow B/J\otimes_A A/I\otimes_A A^{(-1)}$ 
is also injective.
In other words, 
$F_A: B/\fa\rightarrow (B\otimes_A A^{(-1)})/\fa(B\otimes_A A^{(-1)})$ 
is injective.

Let $L\subset \kappa(\fm)^{(-1)}$ 
be a finite extension field of $\kappa(\fm)$.
As $B/\fm B\otimes_{\kappa(\fm)}L$ is CMFI,
\[
F:B/(\fm B+ J)\otimes_{B/\fm B} B/\fm B\otimes_{\kappa(\fm)}L
\rightarrow
B/(\fm B+ J)\otimes_{B/\fm B} (B/\fm B\otimes_{\kappa(\fm)}L)^{(-1)}
\]
is injective.
Taking the inductive limit,
\begin{multline*}
F:B/(\fm B+ J)\otimes_{B/\fm B} B/\fm B\otimes_{\kappa(\fm)}\kappa(\fm)^{(-1)}
\rightarrow\\
B/(\fm B+ J)\otimes_{B/\fm B} (B/\fm B\otimes_{\kappa(\fm)}
\kappa(\fm)^{(-1)})^{(-1)}
\end{multline*}
is injective.
By \cite[Lemma~4.1]{Hashimoto}, {\bf 7}, 
\begin{multline*}
1\otimes
\Psi_1(\kappa(\fm),B/\fm B): 
B/(\fm B+ J)\otimes_{B/\fm B}B/\fm B\otimes_{\kappa(\fm)}\kappa(\fm)^{-1}
\rightarrow 
\\
B/(\fm B+ J)\otimes_{B/\fm B}(B/\fm B)^{-1}
\end{multline*}
is injective.
As $B\otimes_A A^{(-1)}$ and $B^{(-1)}$ are $A^{(-1)}$-flat,
$B/IB\otimes_{A/I}A^{(-1)}/IA^{(-1)}$ and 
$B^{(-1)}/IB^{(-1)}$ are flat modules over the artinian local ring
$A^{(-1)}/IA^{(-1)}$.
Clearly, $y_1,\ldots,y_m$ is a $B/\fm B\otimes_{\kappa(\fm)}
\kappa(\fm)^{(-1)}$-sequence.
Moreover, 
\[
A^{(-1)}/\fm A^{(-1)} \rightarrow  B^{(-1)}/\fm B^{(-1)}
\]
is a flat homomorphism with an $m$-dimensional Cohen--Macaulay closed fiber.
As $A^{(-1)}/\fm A^{(-1)}$ is artinian, 
$B^{(-1)}/\fm B^{(-1)}$ is $m$-dimensional Cohen--Macaulay.
It is easy to see that $B^{(-1)}/(\fm B^{(-1)}+J B^{(-1)})$ is artinian, 
so $y_1,\ldots,y_m$ is a $B^{(-1)}/\fm B^{(-1)}$-sequence.
So by Lemma~\ref{key.thm}, 
\[
1\otimes \Psi_1(A,B):
(B\otimes_A A^{(-1)})/\fa(B\otimes_A A^{(-1)})
\rightarrow
B^{(-1)}/\fa B^{(-1)}
\]
is injective.

Hence $1\otimes F_B: B/\fa \rightarrow B^{(-1)}/\fa B^{(-1)}$ 
is injective.
By Corollary~\ref{CMFI-equiv.thm}, $B$ is CMFI.
\qed

\begin{corollary}
Let $(B,\fn)$ be a noetherian local ring of characteristic $p$,
and $t\in\fn$ be a nonzerodivisor.
If $B/tB$ is CMFI, then $B$ is CMFI.
\end{corollary}

\proof Let $A$ be the localization $\Bbb F_p[T]_{(T)}$ of the 
polynomial ring $\Bbb F_p[T]$ at the maximal ideal $(T)$.
Then the canonical map $A\rightarrow B$ which maps $T$ to $t$ is flat, 
as $t$ is a nonzerodivisor.
Let $L$ be a finite extension field of $\Bbb F_p$.
Then $L$ is a separable extension of $\Bbb F_p$, and hence
$B/tB\otimes_{\Bbb F_p}L$ is \'etale over $B/tB$, and hence is CMFI
by Proposition~\ref{base-change-CMFI.thm}.
Thus $B/tB$ is geometrically CMFI over $A/tA$.
On the other hand, $A$ is regular, and hence is CMFI.
By Proposition~\ref{base-change-CMFI.thm} again, $B$ is CMFI.
\qed

\begin{corollary}\label{CMFI-open.thm}
Let $R$ be an excellent local ring of characteristic $p$, and
$A$ an $R$-algebra essentially of finite type.
Then $\CMFI(A)$ is a Zariski open subset of $\Spec A$.
\end{corollary}

\proof Let $\hat R$ be the completion of $R$, and $\hat A:=\hat R\otimes_R A$.
Note that $\CMFI(\hat A)$ is a Zariski open subset of $\Spec \hat A$
by Lemma~\ref{Gamma-F-injective.thm} and Lemma~\ref{F-finite-open.thm}.

Let $\rho:\Spec \hat A\rightarrow \Spec A$ be the morphism associated with
the base change of the completion.
Then $\rho^{-1}(\NonCMFI(A))=\NonCMFI(\hat A)$ by
Proposition~\ref{base-change-CMFI.thm} and Lemma~\ref{ff-cmfi.thm}.
By Corollary~\ref{CMFI-localization.thm} and 
Lemma~\ref{closed.thm}, $\NonCMFI(A)=\rho(\NonCMFI(\hat A))$ is closed,
as desired.
\qed

\section{Matijevic--Roberts type theorem}

\paragraph 
Let $S$ be a scheme, $G$ an $S$-group scheme, and $X$ a standard
$G$-scheme \cite[(2.18)]{HO2} (that is, $X$ is noetherian and the
second projection $p_2:G\times X\rightarrow X$ is flat of finite type).

\begin{theorem}\label{M-R.thm}
Let $y$ be a point of $X$, and $Y$ the integral closed subscheme of $X$ whose
generic point is $y$.
Let $\eta$ be the generic point of an irreducible component of $Y^*$,
where $Y^*$ is the smallest $G$-stable closed subscheme of $X$ containing
$Y$.
Assume either that the second projection $p_2:G\times X\rightarrow X$ is
smooth, or that $S=\Spec k$ with $k$ a perfect field and $G$ is of finite
type over $S$.
Assume that $\O_{X,\eta}$ is of characteristic $p$.
Then $\O_{X,y}$ is of characteristic $p$.
Moreover, 
\begin{description}
\item[1] If $\O_{X,\eta}$ is $F$-pure, then $\O_{X,y}$ is $F$-pure.
\item[2] If $\O_{X,\eta}$ is excellent and strongly $F$-regular, then
$\O_{X,y}$ is strongly $F$-regular.
\item[3] If $\O_{X,\eta}$ is CMFI, then
$\O_{X,y}$ is CMFI.
\end{description}
\end{theorem}

\proof 
We set $\Cal C$ and $\Cal D$ to be the
class of all noetherian local rings of characteristic $p$, and
$\Bbb P(A,M)$ to be ``always true'' in \cite[Corollary~7.6]{HM}.
Then the conditions {\bf (i)} and {\bf (ii)} there are satisfied,
and by \cite[Corollary~7.6]{HM}, $\O_{X,y}$ is of characteristic $p$.

{\bf 1}
We set $\Cal C$ and $\Cal D$ to be the
class of all $F$-pure noetherian local rings of characteristic $p$, and
$\Bbb P(A,M)$ to be ``always true'' in \cite[Corollary~7.6]{HM}.
Then {\bf (i)} there (the smooth base change) is satisfied by
Lemma~\ref{basic.thm}, {\bf 4} and {\bf 6}.
The condition {\bf (ii)} (the flat descent) holds by Lemma~\ref{basic.thm},
{\bf 5}.
So the assertion follows by \cite[Corollary~7.6]{HM}.

{\bf 2} Set $\Cal C$ to be the class of excellent strongly $F$-regular
noetherian local domains of characteristic $p$, 
and $\Cal D$ to be the class of strongly $F$-regular
noetherian local domains of characteristic $p$.
Then {\bf (i)} and {\bf (ii)} of \cite[Corollary~7.6]{HM}
are satisfied by Lemma~\ref{smooth-base-change.thm} and
Lemma~\ref{pure-sfr.thm}.

{\bf 3} Set $\Cal C=\Cal D$ be the class of CMFI noetherian local
rings of characteristic $p$.
Then {\bf (i)} and {\bf (ii)} of \cite[Corollary~7.6]{HM} are
satisfied by Proposition~\ref{base-change-CMFI.thm} and
Lemma~\ref{ff-cmfi.thm}.
\qed

\begin{corollary}\label{M-R-graded.thm}
Let $p$ be a prime number, and
$A$ a $\Bbb Z^n$-graded noetherian ring of characteristic  $p$.
Let $P$ be a prime ideal of $A$, and $P^*$ be the prime ideal of $A$
generated by the homogeneous elements of $P$.
If $A_{P^*}$ is $F$-pure \(resp.\ excellent strongly $F$-regular,
CMFI\), then
$A_P$ is $F$-pure \(resp.\ strongly $F$-regular, 
CMFI\).
\end{corollary}

\proof Let $S=\Spec \Bbb Z$, $G=\Bbb G_m^n$, and $X=\Spec A$.
If $y=P$, then $\eta$ in Theorem~\ref{M-R.thm} is $P^*$.
The assertion follows immediately by Theorem~\ref{M-R.thm}.
\qed

\begin{corollary}\label{graded-enough.thm}
Let $A$ be a $\Bbb Z^n$-graded noetherian ring of characteristic $p$.
If $A_\fm$ is $F$-pure \(resp.\ excellent strongly $F$-regular, 
CMFI\) for any maximal graded
ideals \(that is, $G$-maximal ideals for $G=\Bbb G_m^n$
\(called *maximal ideal in {\rm\cite{BH}}\)\),
then $A$ is $F$-pure \(resp.\ strongly $F$-regular, Cohen--Macaulay
$F$-injective\).
\end{corollary}

\proof Similar to \cite[Corollary~7.11]{HM}.
\qed

\begin{corollary}\label{graded-deformation.thm}
Let $A=\bigoplus_{n\geq 0}A_n$ be an $\Bbb N$-graded 
noetherian
ring of characteristic $p$.
Let $t\in A_+:=\bigoplus_{n>0}A_n$ be a nonzerodivisor of $A$.
If $A/tA$ is CMFI, then $A$ is 
CMFI.
\end{corollary}

\proof Similar to \cite[Corollary~7.13]{HM}.
\qed

\begin{corollary}\label{filtration-deformation.thm}
Let $A$ be a ring of characteristic $p$, and $(F_n)_{n\geq 0}$ 
a filtration of $A$.
That is, $F_0\subset F_1\subset F_2\subset\cdots\subset A$, 
$1\in F_0$, $F_iF_j\subset F_{i+j}$, and $\bigcup_{n\geq 0}F_n=A$.
Set $R=\bigoplus_{n\geq 0}F_n t^n\subset A[t]$, and $G=R/tR$.
If $G$ is noetherian and Cohen--Macaulay $F$-injective, then
$A$ is also noetherian and Cohen--Macaulay $F$-injective.
\end{corollary}

\proof Similar to \cite[Corollary~7.14]{HM}.
\qed

\end{document}